\documentclass[5p,twocolumn,compress]{elsarticle}
\usepackage{epsfig,graphicx,color}
\usepackage{amsmath, amsthm, amsfonts,amssymb}
\usepackage{dcolumn}
\usepackage{empheq}
\usepackage{tikz}
\usepackage{ulem}
\usepackage{colortbl}
\usepackage{comment}
\usepackage{graphicx}
\usepackage{txfonts}
\usepackage{float}
\newenvironment{acknowledgements}
{\section*{Acknowledgements}}
{}

\definecolor{Green}{rgb}{0.0, 0.5, 0.0} 
\journal{Chaos, Solitons \& Fractals}
\begin{document}
\begin{frontmatter}

\title{Hysteresis in a Generalized Kuramoto Model \\ with a Simplified Realistic Coupling Function and Inhomogeneous Coupling Strengths}
\author[inst1]{Jae Hyung Woo}
\affiliation[inst1]{ 
organization={Department of Psychological and Brain Sciences, Dartmouth College}, 
city={Hanover},
postcode={03766},
state={NH},
country={USA}
}
\author[inst2]{Hae Seong Lee}
\affiliation[inst2]{ 
organization={Department of Physics, Sungkyunkwan University}, 
city={Suwon},
postcode={16419}, 
country={Republic of Korea}
}

\author[inst3,inst4]{Joon-Young Moon\corref{corauth}}
\ead{joon.young.moon@gmail.com}
\affiliation[inst3]{
	organization={Center for Neuroscience Imaging Research, Institute for Basic Science (IBS)}, 
city={Suwon},
postcode={16419}, 
country={Republic of Korea}
}
\affiliation[inst4]{
	organization={Sungkyunkwan University}, 
city={Suwon},
postcode={16419}, 
country={Republic of Korea}
}

\author[inst5]{Tae-Wook Ko\corref{corauth}}
\ead{twko@nims.re.kr}
\affiliation[inst5]{
organization={National Institute for Mathematical Sciences}, 
city={Daejeon},
postcode={34047},
country={Republic of Korea}
}

\cortext[corauth]{Corresponding authors.}

\date{\today}

\begin{abstract}
We investigate hysteresis in a generalized Kuramoto model with identical oscillators, focusing on coupling strength inhomogeneity, which results in oscillators being coupled to others with varying strength, and a simplified, more realistic coupling function.  With the more realistic coupling function and the coupling strength inhomogeneity, each oscillator acquires an effective intrinsic frequency proportional to its individual coupling strength. This is analogous to the positive coupling strength-frequency correlation introduced explicitly or implicitly in some previous models with nonidentical oscillators that show explosive synchronization and hysteresis. Through numerical simulations
and analysis using truncated Gaussian, uniform, and truncated power-law coupling strength distributions, we observe that the
system can exhibit abrupt phase transitions and hysteresis. The distribution of coupling strengths significantly affects the hysteresis
regions within the parameter space of the coupling function. 
Additionally, numerical simulations of models with weighted networks including a brain network confirm the existence of hysteresis due to the realistic coupling function and coupling strength inhomogeneity, suggesting the broad applicability of our findings to complex real-world systems.
\end{abstract}

\end{frontmatter}

\section{\label{sec:Introduction}Introduction}
Hysteresis, a phenomenon where the state of a system depends on its history \cite{visintin1994,chow2023}, has been widely observed in various systems, including physical \cite{bertotti1998}, biological \cite{ferrell2001,tyson2003,chen2000,cross2002,noori2014} and neural systems \cite{izhikevich2007,kleinschmidt2002,hkim2018,sepulveda2018}.
In physics, it is fundamental for comprehending the properties of magnetic materials and phase transitions \cite{bertotti1998}. Similarly, in the biological realm, hysteresis plays a crucial role in the functions of regulatory and signaling systems of cells \cite{ferrell2001,tyson2003,chen2000,cross2002,noori2014}.  In neuroscience, hysteresis is important as it explains how neurons, the brain, and their functions are affected by the historical activity, particularly in responses to stimuli or anesthesia \cite{izhikevich2007,kleinschmidt2002,hkim2018,sepulveda2018}. 

Studies focusing on coupled oscillator systems, which serve as  models for rhythmic activities in complex systems such as the brain \cite{pikovsky2001,strogatz2003,kuramoto1984,strogatz2000,acebron2005,ermentrout2001}, have shown several settings for hysteresis \cite{schuster1989,skim1997,yeung1999,twko2004,twko2007,boccaletti2016}. 
In the first type of setting, time delay in interactions between oscillators induces bistability or multistability and hysteresis~\cite{schuster1989,skim1997,yeung1999,twko2004,twko2007}.  
In the second type of setting, a positive correlation between the total coupling strengths ($K_i$) or the number ($k_i$) of connections to individual oscillators and the absolute values of the intrinsic frequencies of the oscillators is essential for the manifestation of an abrupt phase transition to synchronization, called explosive synchronization (ES), and hysteresis \cite{gomez2011,xzhang2013,xzhang2014,vlasov2015,xiao2017,kundu2017,kundu2019}. Suppression of the merging of small synchronized clusters into a giant cluster was suggested as the reason for the abrupt phase transition and thus the hysteresis \cite{xzhang2014}.
In the third type of setting where the synchronization route via cluster formation is avoided by weakening or removing the links between oscillators with similar frequencies, ES and hysteresis are observed \cite{leyva2013a,leyva2013b,lzhu2013}. Specifically, in Ref. \cite{leyva2013b}, the coupling strengths between connected oscillators $i$ and $j$ are given by an increasing function of the intrinsic frequency mismatch $|\omega_i - \omega_j|$.
The weighting leads to stronger coupling strengths for oscillators with frequencies lying in less populated regions of the frequency distribution—typically the lower and higher frequency regions for symmetric distributions, and the higher frequency region for decreasing distributions \cite{leyva2013b}. This spontaneously obtained property corresponds to the imposed correlation between $K_i$ and $|\omega_i|$ of the second type of setting. 
In the fourth type of setting, the system partly or fully employs a dynamic coupling strength that is proportional to the extent of synchronization between oscillators. Thus, the synchronizing force is of an order greater than one with respect to the extent of synchronization \cite{filatrella2007,xzhang2015}. As a result, once synchronization is achieved, the coupling strength must be significantly lowered to disrupt this synchronized state. This type of system shows ES and hysteresis, and it was emphasized in Ref.~\cite{xzhang2015} that the necessary condition for the existence of ES (thus hysteresis) is the presence of a suppressive rule. This rule must be able to prevent the formation of a giant synchronization cluster, regardless of its specific nature.

Interestingly, in different scenarios involving variations in coupling topology and coupling functions, ranges of parameters exhibiting various forms of bistability have been observed in coupled oscillator systems, although hysteresis was not directly examined in these cases. In nonlocally coupled identical oscillators with a phase-shifted coupling function, bistability between a fully locked state and a partially locked state, known as a chimera state, has been observed~\cite{kuramoto2002, abrams2004, abrams2006, abrams2008}. Similarly, in coupled identical oscillators with coupling strength inhomogeneity and a more realistic coupling function, two types of bistability have been reported: between a partially (or fully) locked state and an incoherent state~ \cite{twko2008a}, and between a fully locked state and a partially locked state~\cite{twko2008b}.  
It is highly probable that hysteresis exists in these systems since bistability typically implies hysteresis. 
We need to study hysteresis in such simple systems of coupled oscillators without any special factors.

To address this issue, in this paper we investigate the presence of hysteresis in the most fundamental system of coupled oscillators characterized solely by the essential elements: a coupling function and the inherent inhomogeneity in coupling strengths.   
The coupling function we use is a more realistic one simplified from realistic coupling functions.
Our aim is to elucidate how these two fundamental factors impact hysteresis behavior. We deliberately exclude any additional elements such as time delays, coupling strength-frequency correlations, intended weighting of links, or synchronization-dependent coupling strength explained above. By doing so, we examine hysteresis in the simplest possible configuration of coupled oscillators, providing insights into the core dynamics of these systems.
Utilizing a mean-field model of coupled identical oscillators based on the Kuramoto framework, we demonstrate that hysteresis manifests when strength inhomogeneity is present alongside with a more realistic coupling function.
We find that both the distribution of coupling strengths and the specific nature of the coupling functions critically influence the emergence and characteristics of hysteresis within such systems.
Numerical simulations with a network version of the model, specifically one using a brain network, generalize our findings to the cases of complex real-world systems.

\section{\label{sec:Model}Model}
Here, we analyze a generalized mean-field Kuramoto model with a simplified realistic coupling function and inhomogeneous coupling strengths introduced and studied in Refs. \cite{twko2008a,twko2008b,jkim2019,laing2009} to elucidate the effects of coupling strength inhomogeneity and the coupling function on hysteresis in systems of coupled oscillators. 
\begin{eqnarray}
	\frac{d \theta_i}{dt} &=& \omega + \frac{K_i}{N}\sum_{j=1}^N \big [c_0+\sin(\theta_j - \theta_i-\beta)\big],  
\label{eq:model} \\
	&&i=1,2,...,N, ~~c_0 \geq 0, ~~\beta \in [0,\pi/2), \nonumber
\end{eqnarray}
where $\theta_i(t)$ is the phase of oscillator $i$ at time $t$, $\omega$ is the intrinsic frequency, $N$ is the total number of oscillators. 
$K_i$ denotes the total coupling strength to oscillator $i$ from all other oscillators, and is a positive number randomly drawn from a given probability distribution $g(K)$.
This model encapsulates a phase-reduced system of weakly coupled, identical limit-cycle oscillators, characterized by their (in)homogeneous coupling strengths. The coupling function $H(\theta) = c_0 + \sin(\theta - \beta)$ is derived as a first-order approximation of general coupling functions via the phase reduction method~\cite{kuramoto1984,strogatz2000,acebron2005,ermentrout2001,twko2008a,mhpark1996}. The coefficient of the second term of $H$ is set to $1$ without loss of generality. 
This coupling function satisfies the condition $H'(0)>0$, with which the system shows locally stable in-phase synchronous behavior in the absence of coupling strength inhomogeneity \cite{ermentrout1992}.
Notably, the constant $c_0$, which has traditionally been omitted in the Kuramoto model and its variants, is retained in our model. Specifically, the Kuramoto model employs a coupling function $H_{K}(\theta) = \sin\theta$~\cite{kuramoto1984,strogatz2000,acebron2005}, while the Sakaguchi-Kuramoto model uses $H_{KS}(\theta) = \sin(\theta - \beta)$~\cite{sakaguchi1986}. The absence of inhomogeneity in the coupling strength $K_i$ renders $c_0$ non-influential to the dynamics, thus justifying its exclusion in previous models. However, in our framework, both constants $c_0$ and $\beta$ play crucial roles in shaping the coupling function $H$, underscoring the importance of considering (in)homogeneous coupling strengths.

Previous studies with the model showed that the model can exhibit various states including uniformly incoherent state and partially locked states \cite{twko2008a,twko2008b,jkim2019}. With uniform coupling strength distributions, bistability between a partially (fully) locked state and a uniformly incoherent state was observed in this model~\cite{twko2008a}. 

Let us consider a situation where we change a parameter of the full model of weakly coupled limit-cycle oscillators, such as the input to the oscillators, excluding the coupling strengths. This parameter change causes changes in the parameters $\omega$, $c_0$, and $\beta$ of the model (Eq.~(\ref{eq:model})) through the phase-reduction. 
Therefore, choosing such a parameter of the full model as the control parameter for the hysteresis procedure is equivalent to using the parameters $(\omega, c_0, \beta)$ collectively as the control parameter along a specific path in the parameter space of the reduced model. 
Even though $c_0$ of our model cannot typically be adjusted independently through the change of the parameters of the full model, we treat $c_0$ as an independent control parameter for hysteresis procedures, while holding $\beta$ constant. This approach can provide insights into hysteresis behaviors along various paths in the parameter space of the reduced model.

Note that in the previous studies of explosive synchronization and hysteresis the global coupling strength affecting all the oscillators is used as the control parameter for hysteresis procedures \cite{gomez2011,xzhang2013,xzhang2014,vlasov2015,xiao2017,kundu2017,kundu2019,leyva2013a,leyva2013b,lzhu2013,filatrella2007,xzhang2015}. Here, with identical oscillators, introducing such global coupling strength has no qualitative effect. 

We calculate the order parameter $r(t)$ to measure the extent of synchronization and measure its time average $R$ during the hysteresis procedures. They are defined as follows: 
\begin{eqnarray}
\label{eq:order_parameter}
	r(t) e^{i\Theta(t)} &=& \frac{1}{N}\sum_{j=1}^Ne^{i\theta_j(t)}, \\
	R &=& \langle r(t) \rangle_t,
\end{eqnarray}
where $r(t)$ and $\Theta(t)$ are the order parameter and the corresponding phase, respectively. The notation $\langle \cdot\rangle_t$ stands for the time average. The value of $r(t)$ ranges from 0 (complete incoherence) to 1 (full synchronization).

To better understand our model, we apply the following transformation $\theta_i \rightarrow \omega t + \phi_i$, $t \rightarrow \tau /c_0$ with $c_0 > 0$, and get the following from Eq. (\ref{eq:model}) 
\begin{equation}
\label{eq:model_tr}
	\frac{d \phi_i}{d\tau} =K_i+\frac{1}{c_0}\frac{K_i}{N}\sum_{j=1}^N \sin(\phi_j - \phi_i-\beta).
\end{equation}
This transformed model has the form of coupled nonidentical oscillators with intrinsic frequencies $\omega_i = K_i$, the global coupling strength $1/c_0$ affecting all oscillators, and the local coupling strength $K_i$. The natural setting, incorporating a more realistic coupling function and coupling strength inhomogeneity, results in a positive correlation between the effective intrinsic frequencies $\omega_i$ and the coupling strength $K_i$. This correlation is similar to that found in previous models with nonidentical oscillators that exhibit explosive synchronization and hysteresis \cite{gomez2011,xzhang2013,xzhang2014,vlasov2015,xiao2017,kundu2017,kundu2019,leyva2013a,leyva2013b,lzhu2013}.
While this transformed model is different from the Sakaguchi-Kuramoto model with positive and negative intrinsic frequencies $\omega_i$ and $K_i=|\omega_i|$ of Ref.~\cite{xiao2017}, it bears a particular resemblance to the Sakaguchi-Kuramoto model on scale-free networks with degree-frequency correlation~\cite{kundu2017,kundu2019}, where degrees and frequencies are positive, as in our model. 
From the resemblance, we can expect explosive synchronization-like abrupt phase transitions and hysteresis in our model with identical oscillators.

However, a crucial difference emerges when considering $1/c_0$ as the effective coupling strength. Unlike in typical models of nonidentical oscillators where the system transitions from an incoherent state to synchronization \cite{gomez2011,xzhang2013,xzhang2014,vlasov2015,xiao2017,kundu2017,kundu2019,leyva2013a,leyva2013b,lzhu2013}, our model exhibits a reversed hysteresis phenomenon. Specifically, as $c_0$ increases (and thus $1/c_0$ decreases), the system transitions from a coherent state to an incoherent state. This reversed behavior distinguishes our model while still maintaining the possibility of explosive synchronization-like abrupt phase transitions and hysteresis.

\section{\label{sec:Analysis}Analysis}
Since hysteresis is a phenomenon of history-dependent selection of the state of a system, we need to know the stable states of the model.
To simplify our investigation, we focus on stationary states, where $r(t)$ and $\Theta(t)$ remain constant over time, and find out these states and their stability through following analyses.     

\subsection{Self-consistency analysis for finding stationary states}
To find stationary states, we conduct a self-consistency analysis, which has been utilized to identify such states in the Kuramoto model and its variants \cite{kuramoto1984,strogatz2000,acebron2005,kuramoto2002,abrams2004,abrams2006,abrams2008,twko2008b,jkim2019}. In this analysis, we calculate the order parameter self-consistently by aggregating contributions from both locked and drifting oscillators. The behavior of these oscillators is, in turn, determined by the order parameter itself.

For the model of Eq.~(\ref{eq:model}), self-consistency approach gives the following equation obtained in Refs. \cite{twko2008b,jkim2019}.
\begin{eqnarray}
\label{eq:self_consistency}
	r^2 {\rm e}^{i\beta} &=& i\int_{\mathcal{D}_{tot}} dK \: \frac{g(K)Z(K)}{K}  
\nonumber \\
	&&+\int_{\mathcal{D}_l} dK\: \frac{g(K)\sqrt{K^2 r^2 -Z(K)^2}}{K} \\
	&& -i\int_{\mathcal{D}_d} dK\: \frac{g(K){\rm sign}\big(Z(K)\big)\sqrt{Z(K)^2-K^2 r^2}}{K},\nonumber
\end{eqnarray}
where $\mathcal{D}_{\text{tot}}$ represents the total range of the coupling strength $K$, encompassing all oscillators. The subset $\mathcal{D}_{l}$ is defined as the range(s) of $K$ for which oscillators are phase-locked at a common frequency $\Omega$, indicating synchronous behavior. Conversely, $\mathcal{D}_{d}$ denotes the range(s) of $K$ corresponding to oscillators that drift monotonically without achieving phase-locking. We introduce $\Delta \equiv \omega - \Omega$ as the difference between the intrinsic frequency $\omega$ and the common frequency $\Omega$, and define $Z(K) \equiv \Delta + Kc_0$. 
Accordingly, $\mathcal{D}_{l}$ is specified by $\{K_i : K_i \:r > |\Delta + K_i c_0|\,\}$, capturing oscillators whose coupling strength facilitates phase-locking, whereas $\mathcal{D}_{d}$ is determined by $\{K_i : K_i \:r < |\Delta + K_i c_0|\,\}$, identifying those that drift.

Note that self-consistency analysis does not determine the stability of the obtained states.

We numerically search for states represented by $(r,\Delta)$ that satisfy the above equation for given $g(K)$, $c_0$, and $\beta$. In this search, we explicitly consider the possibility of multiple solutions, which is crucial for understanding hysteresis in the system. We subdivide the range $[0,1]$ for $r$ into subintervals to ensure a thorough search for all possible solutions.

\subsection{Stability analysis of the stationary states}
Regarding the stability of a uniformly incoherent state of the model Eq.~(\ref{eq:model}) with $H(\theta)=c_0+\sin(\theta-\beta)$, in Refs.~\cite{twko2008a,twko2008b}, it was shown through linear stability analysis using population density equation that the eigenvalues $\lambda = \mu-i\nu-i\omega$ determining the stability of the state satisfy the equations:
\begin{subequations}
\begin{eqnarray}
	\cos \beta &=& \frac{1}{2} \int_0^\infty dK \: \frac{\mu K g(K)}{\mu^2+(Kc_0-\nu)^2}, 
\label{eq:stability_a}
\\
	\sin \beta &=& \frac{1}{2} \int_0^\infty dK \: \frac{K g(K)(Kc_0-\nu)}{\mu^2+(Kc_0-\nu)^2}.
\label{eq:stability_b}
\end{eqnarray}
\end{subequations}

In the limit of $\mu\rightarrow 0^+$, we can identify a critical threshold for $c_0$ below which the uniformly incoherent state is unstable. The state is linearly neutrally stable above the threshold $c_0$ \cite{twko2008a}. 
In this limit, Eqs.~(\ref{eq:stability_a}) and (\ref{eq:stability_b}) become the following equations
\begin{subequations}
\begin{eqnarray}
\cos \beta &=&  \frac{\pi}{2 c_0} \bigg (\frac{\nu}{c_0} \bigg ) g\bigg (\frac{\nu}{c_0} \bigg ), 
\label{eq:uniform_state_stability_a}
\\
\sin \beta &=& \lim_{\mu\rightarrow 0^+}\frac{1}{2c_0} \int_{K_a}^{K_b} dK \: \frac{K g(K)(K-\frac{\nu}{c_0})}{(\frac{\mu}{c_0})^2+(K-\frac{\nu}{c_0})^2} \nonumber \\
&&{\rm with} ~0<K_a < \frac{\nu}{c_0} < K_b, 
\label{eq:uniform_state_stability_b}
\end{eqnarray}
\end{subequations}
where $K_a$ and $K_b$ are the lower bound and upper bound for $K$, respectively. Let us call the critical threshold $c_0$ as $c_0^*$.
With a given coupling strength distribution $g(K)$, we numerically obtain $c_0^*$ values as a function of $\beta$. 

For the stability of partially or fully locked states of the system, we use the Ott-Antonsen (OA) approach~\cite{ott2008,ott2009} as in the Refs.~\cite{laing2009,wzou2020,cxu2021} (see \ref{appendix:A} for the details). 
We note that this approach reproduces the previously obtained self-consistency equation (Eq. (\ref{eq:self_consistency})) for the stationary states. Additionally, it recovers Eqs.~(\ref{eq:stability_a}) and (\ref{eq:stability_b}) for the stability of a uniformly incoherent state.
The partially (fully) locked state is unstable if there is an eigenvalue $\lambda$ satisfying the following characteristic equation
\begin{eqnarray}
	(1-J_{11}(\lambda))(1-J_{22}(\lambda)) - J_{12}(\lambda)J_{21}(\lambda) = 0,~~
\label{eq:char_eq_1}
\end{eqnarray}
where
\begin{eqnarray}
	J_{11}(\lambda) &=& \int_0^\infty dK \: \frac{q_{11}(K)g(K)}{\lambda-m_{11}(K)}, \\
	J_{12}(\lambda) &=& \int_0^\infty dK \: \frac{q_{12}(K)g(K)}{\lambda-m_{11}(K)}, \\
	J_{21}(\lambda) &=& \int_0^\infty dK \: \frac{q_{12}^*(K)g(K)}{\lambda-m_{11}^*(K)}, \\
	J_{22}(\lambda) &=& \int_0^\infty dK \: \frac{q_{11}^*(K)g(K)}{\lambda-m_{11}^*(K)}, 
\end{eqnarray}
and the condition $\mathrm{Re}(\lambda) > 0$. We numerically check the existence of such an eigenvalue and determine the stability. If such an eigenvalue is not found, the partially locked state and the fully locked state are linearly neutrally stable and linearly stable, respectively. Other details, including the definitions of $m_{11}(K)$, $q_{11}(K)$ and $q_{12}(K)$, are in \ref{appendix:A}. 

\section{\label{sec:Simulation}Numerical simulation and analysis results for different distributions of coupling strengths}
%=============================
%Fig. 1
\begin{figure*}[bt!]
    \centering
    \includegraphics[width=12cm]{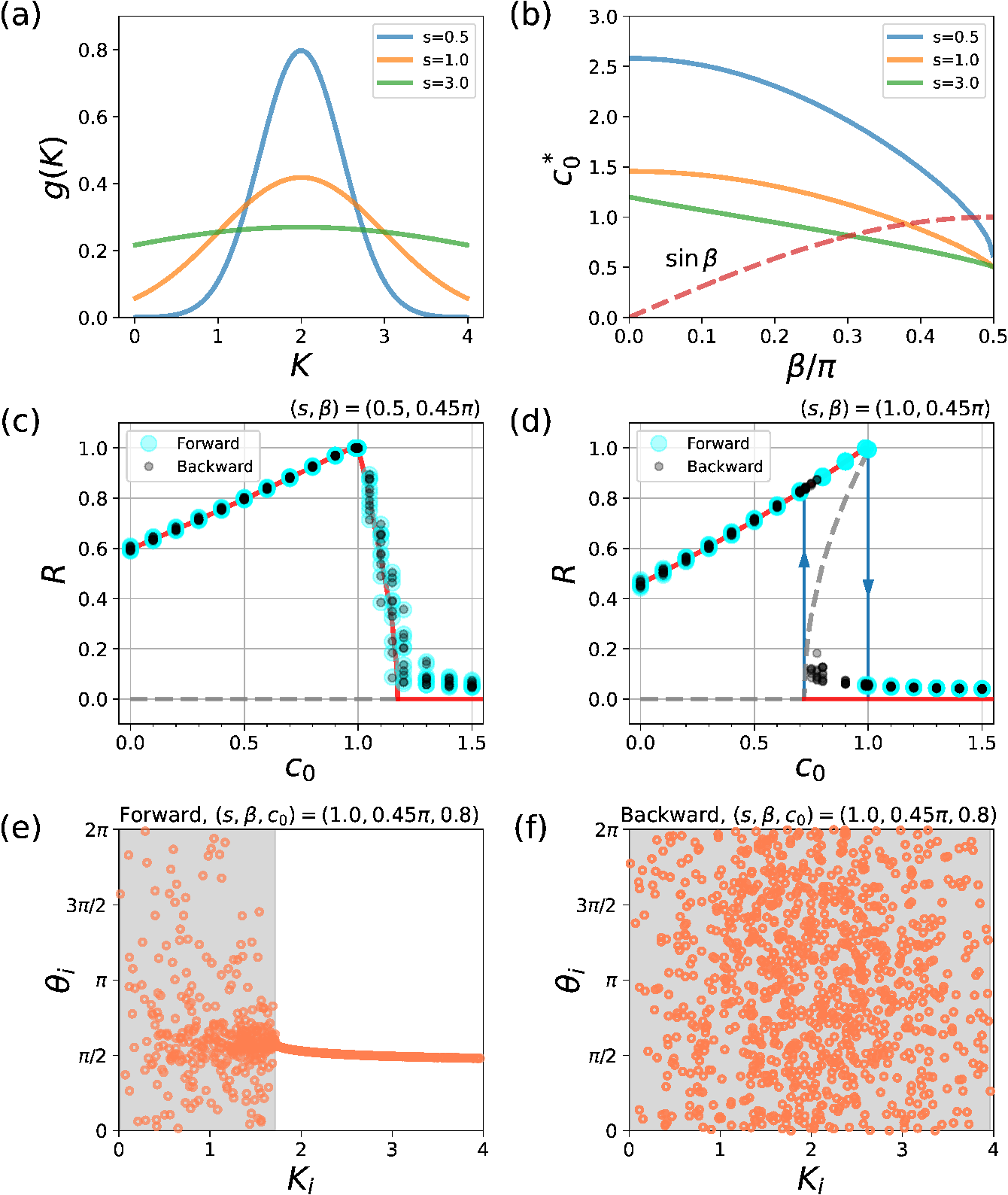}
	\caption{Hysteresis with truncated Gaussian distributions for coupling strength $K_i$. (a) Truncated Gaussian distributions of Eq. (\ref{eq:tr_gauss_dist}) for several values of $s$. (b) Curves of critical $c_0$ ($c_0^*$) for the stability of a uniformly incoherent state as a function of $\beta$ with the distributions of (a). The curves are obtained numerically from Eqs. (\ref{eq:uniform_state_stability_a}) and (\ref{eq:uniform_state_stability_b}). $\sin\beta$ is plotted for the reference. (c) A case with no hysteresis ($s=0.5$, $\beta=0.45\pi$) : $R$ values for the forward path (cyan circles) and for the backward path (gray circles) from 10 hysteresis procedures with different $\{K_i\}$ and initial conditions. See the text for more details. (d) A case with hysteresis ($s=1.0$, $\beta=0.45\pi$). $R$ vs $c_0$ as in (c). 
    In (c) and (d), the red solid lines and the gray dashed lines represent $r$ values obtained numerically from Eq.~(\ref{eq:self_consistency}), corresponding to at least linearly neutrally stable and unstable states, respectively. In (d), the lines with arrows in the middle denote the theoretically obtained abrupt transitions along the hysteresis procedure, considering the stability analysis results. For the case of (d), the system shows two different states with the same parameter values ($s=1.0$, $\beta=0.45\pi$, $c_0=0.8$) along the hysteresis procedure: (e) Forward path: a partially locked state with $R \approx 0.886$ and (f) Backward path: an incoherent state with $R \approx 0.070$.}
\label{fig:01_trg}
\end{figure*}
%=============================

%=============================
%Fig. 2
\begin{figure*}[tb!]
    \centering
    \includegraphics[width=17cm]{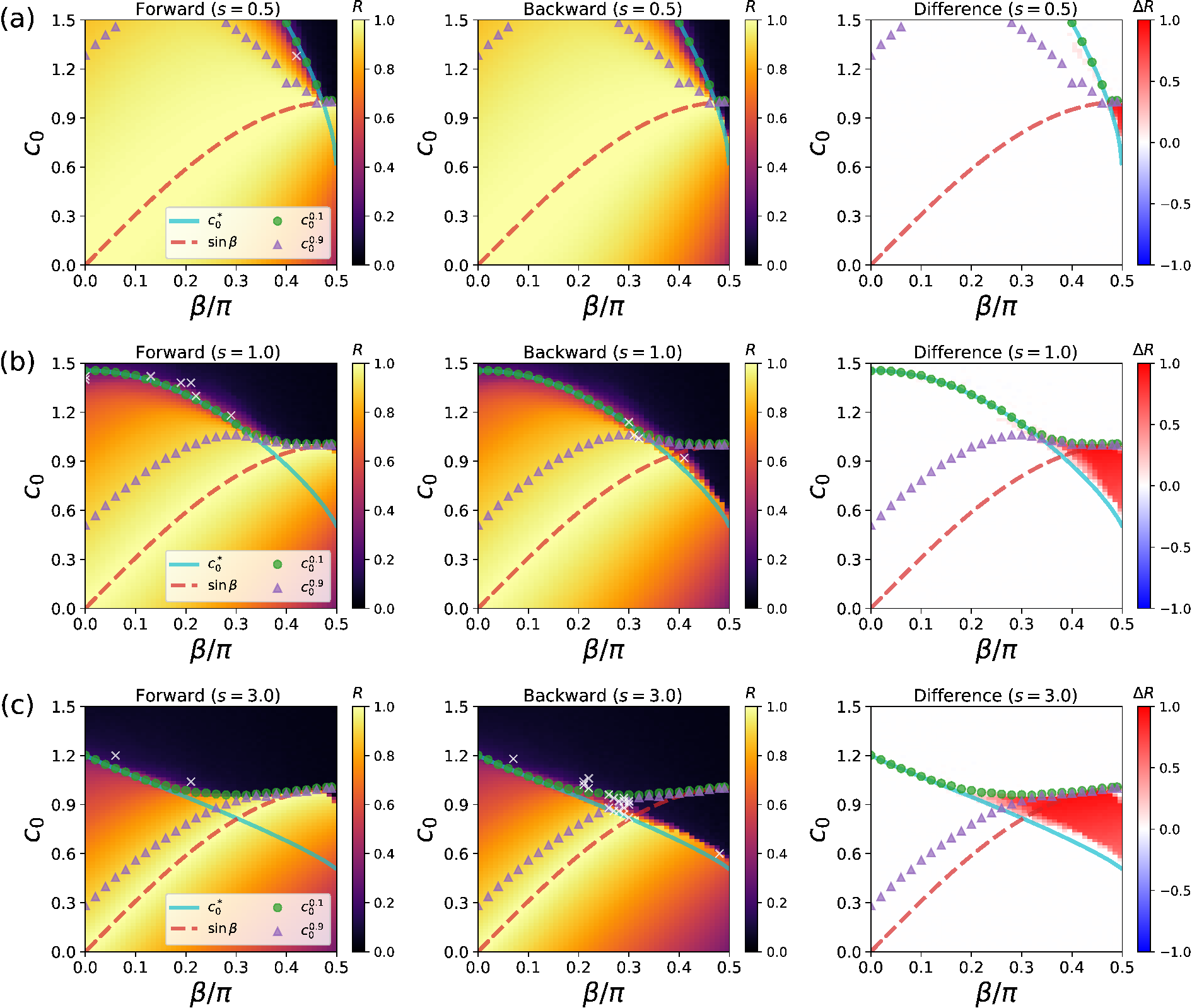}
	\caption{Hysteresis search in $\beta$-$c_0$ plane using hysteresis procedures with truncated Gaussian distributions for every value of $\beta \in [0,~0.49\pi]$ separated by $\Delta \beta=0.01\pi$.  (a) $s=0.5$: (Forward) $R_\text{forward}(\beta,c_0)$, time-averaged order parameter for the forward path, averaged over 5 hysteresis procedures with different $\{K_i\}$ and initial conditions; (Backward) $R_\text{backward}(\beta,c_0)$, $R$ values for the backward path, averaged over the 5 hysteresis procedures; (Difference) $\Delta R(\beta,c_0) \equiv R_\text{forward}(\beta,c_0) - R_\text{backward}(\beta,c_0)$. (b) $s=1.0$. (c) $s=3.0$. $c_0^{0.9}$ and $c_0^{0.1}$ denote the boundaries explained in detail in Subsec.~\ref{subsec:tr_gaussian_dist}. White X marks indicate data points where the standard deviation(std) of the order parameter time series $r(t)$ exceeds $0.1$ for at least one simulation out $5$. See main text for other details.  
	}
\label{fig:02_trg}
\end{figure*}
%=============================

We now simulate and analyze the model Eq. (\ref{eq:model}) for three different coupling strength distributions: truncated Gaussian, uniform, and truncated power-law distributions. These distributions are chosen to represent a range of common scenarios in coupled oscillator systems. We set the total number of oscillators to $N=1000$ and the intrinsic frequency to $\omega=\pi$. Other choices of $\omega$ and $N$, provided $N\geq 1000$, do not change the results qualitatively. For the simulations, we employ the fourth-order Runge-Kutta method with a time step of $\Delta t=0.01$. 
We typically simulate the system for 1000 time units, but this duration may be shorter or longer depending on the specific behavior observed.
In rare cases where values of the order parameter continue to oscillate, we extend simulations up to 15000 time units. 

\begin{figure*}[ht!]
    \centering
    \includegraphics[width=12cm]{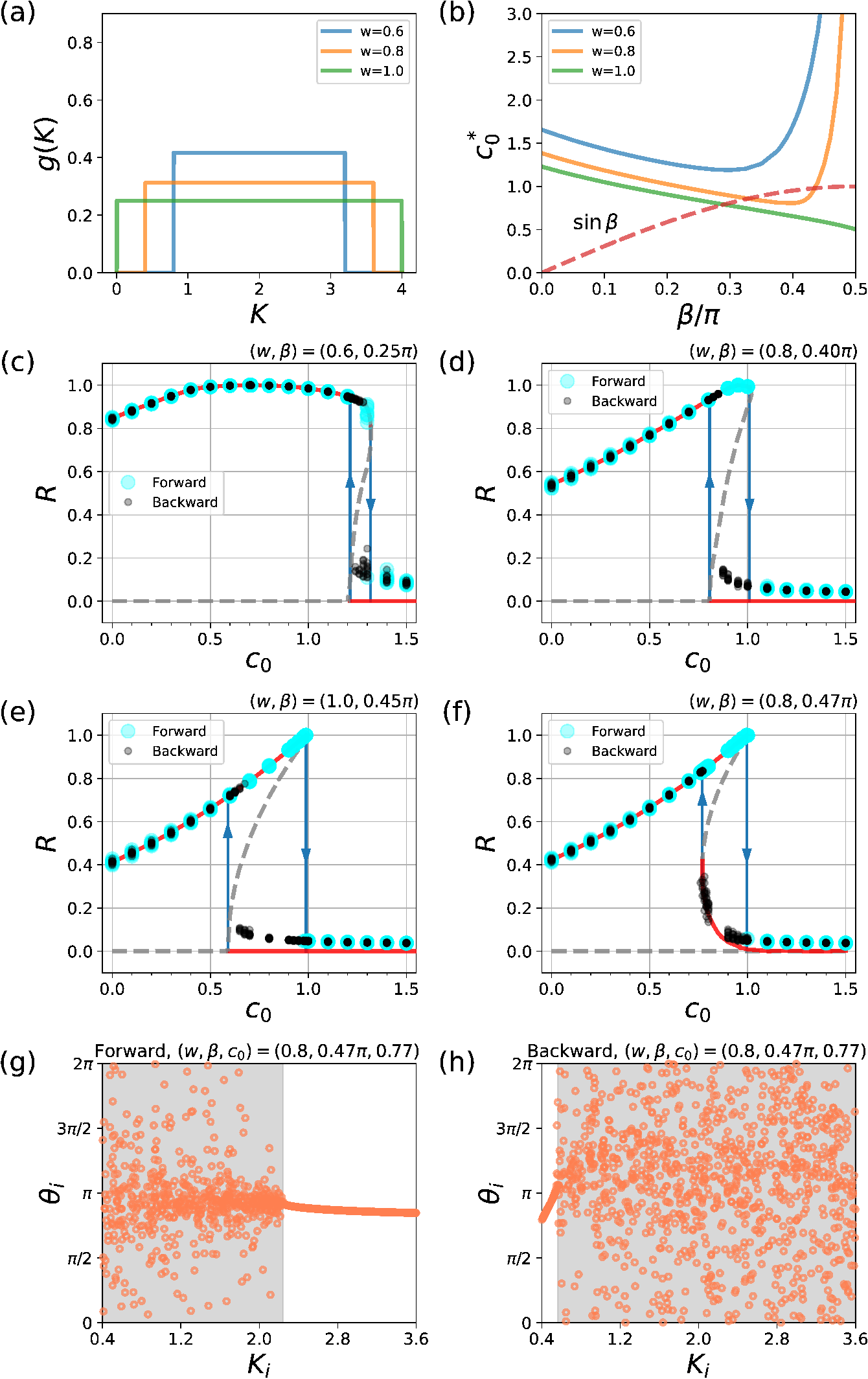}
	\caption{Hysteresis with uniform distributions for coupling strength $K_i$. (a) Uniform distributions of Eq. (\ref{eq:uniform_dist}) for several values of $w$. (b) Curves of critical $c_0$($c_0^*$) as a function of $\beta$ with the distributions of (a).
	$R$ vs $c_0$ for the cases with hysteresis: (c) $w=0.6$, $\beta=0.25\pi$, (d) $w=0.8$, $\beta=0.40\pi$, (e) $w=1.0$, $\beta=0.45\pi$ and (f) $w=0.8$, $\beta=0.47\pi$. 
	For the case of (f), the system shows two different states with the same parameter values ($w=0.8$, $\beta=0.47\pi$, $c_0=0.77$) along the hysteresis procedure: (g) Forward path: a partially locked state with $R\approx 0.837$ and (h) Backward path: a partially locked state with $R \approx 0.296$. Other details are as in Fig \ref{fig:01_trg}.
	}
    \label{fig:01_uni}
\end{figure*}

\begin{figure*}[tb!]
    \centering
    \includegraphics[width=12cm]{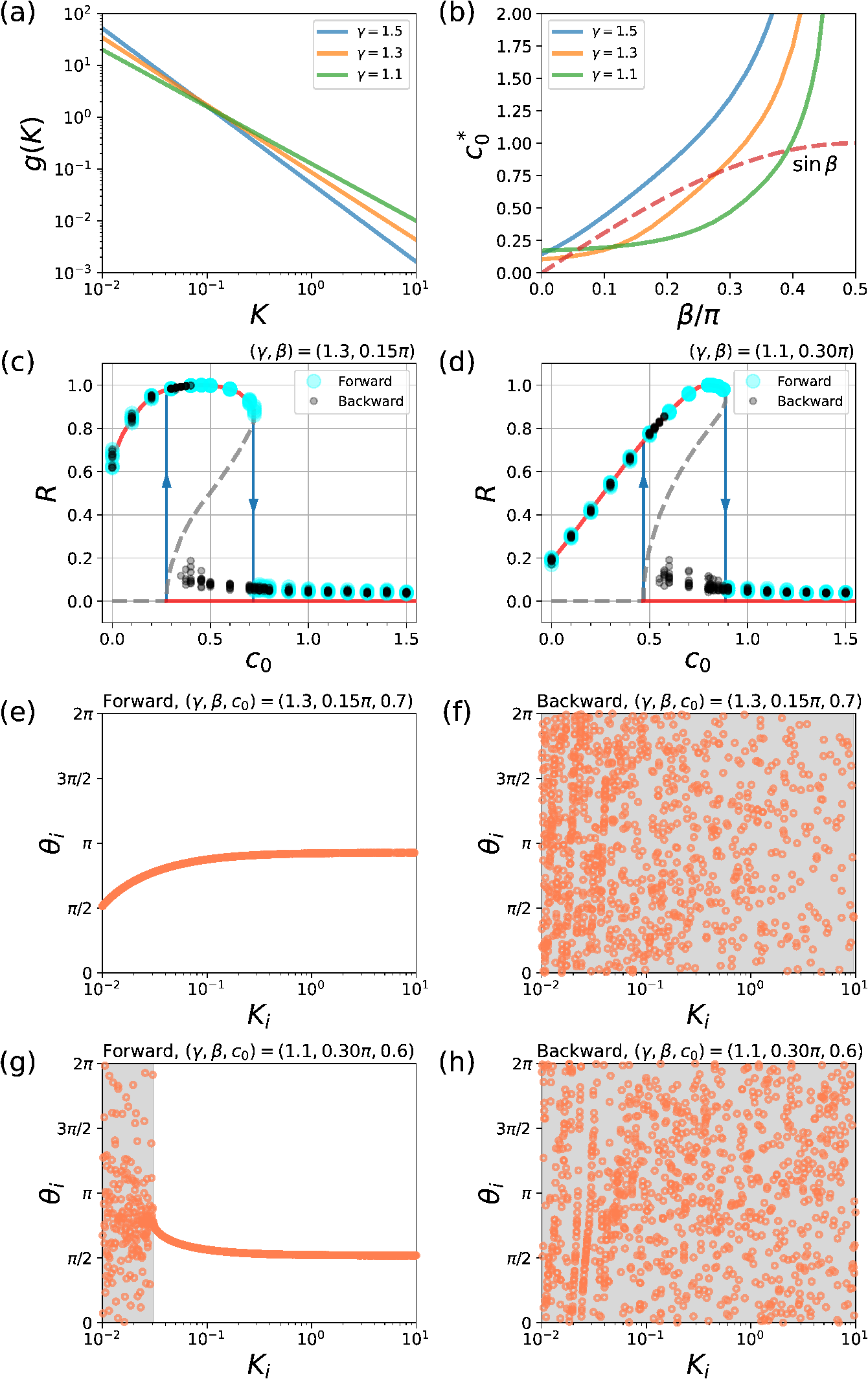}
	\caption{Hysteresis with truncated power-law distributions for coupling strength $K_i$. (a) Truncated power-law distributions of Eq. (\ref{eq:tr_power_law_dist}) for several values of $\gamma$. (b) Curves of critical $c_0$($c_0^*$) as a function of $\beta$ with the distributions of (a). $R$ vs $c_0$ for the cases with hysteresis: (c) $\gamma=1.3$, $\beta=0.15\pi$ and (d) $\gamma=1.1$, $\beta=0.30\pi$.  
	For the case of (c), the system shows two different states with the same parameter values ($\gamma=1.3$, $\beta=0.15\pi$, $c_0=0.7$) along the hysteresis procedure: (e) Forward path: a fully locked state with $R \approx 0.931$ and (f) Backward path: an incoherent state with $R \approx 0.058$.
	For the case of (d), the system shows two different states with the same parameter values ($\gamma=1.1$, $\beta=0.30\pi$, $c_0=0.6$) along the hysteresis procedure: (g) Forward path: a partially locked state with $R \approx 0.874$ and (h) Backward path: an incoherent state with $R \approx 0.087$. Other details are as in Fig \ref{fig:01_trg}.
	}
    \label{fig:01_trpl}
\end{figure*}

\begin{figure*}[bt!]
    \centering
    \includegraphics[width=17cm]{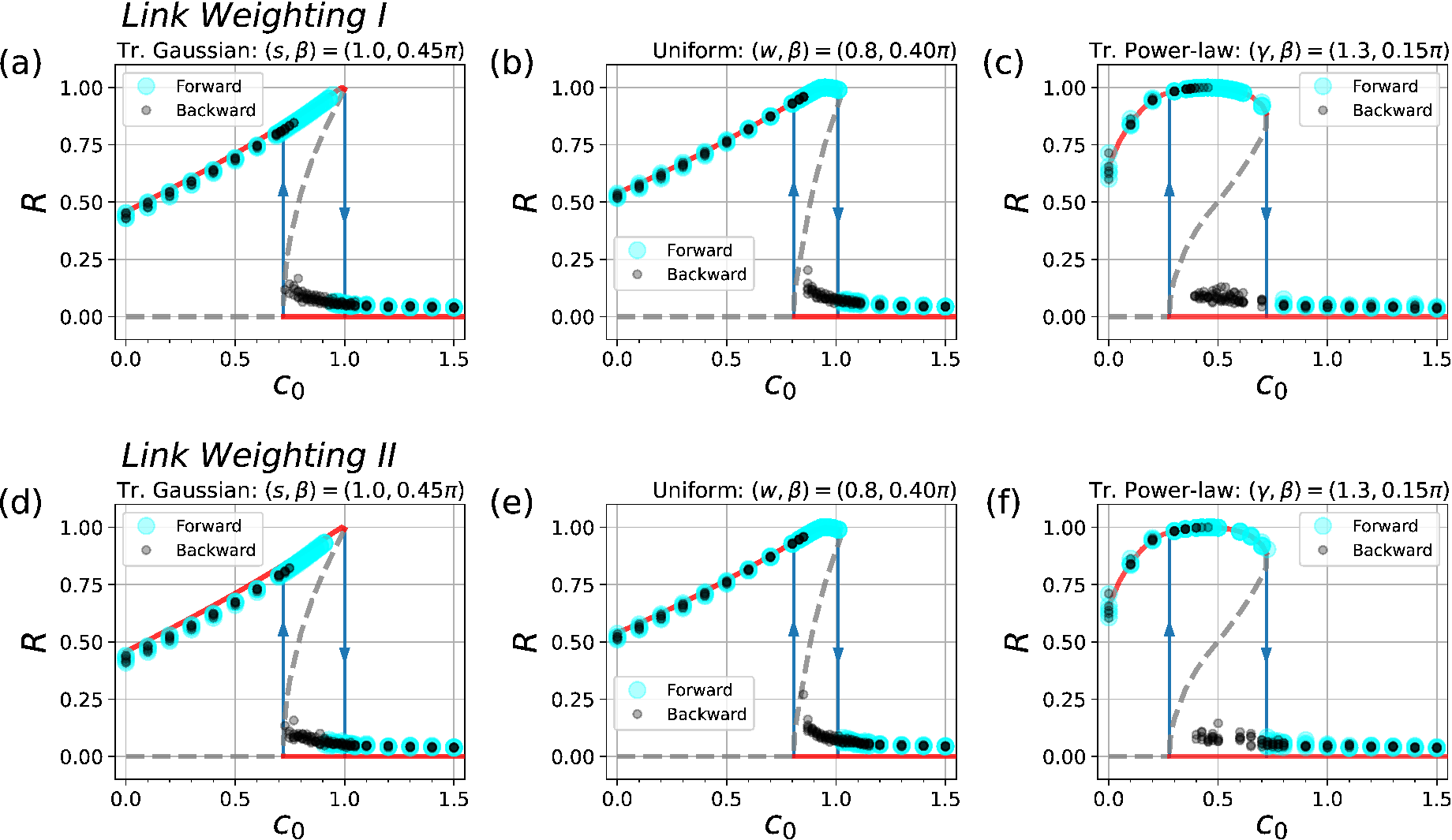}
	\caption{Hysteresis in models with Erd{\H o}s-R{\'e}nyi networks weighted by two different link weighting schemes: Link weighting I (Eq.~(\ref{eq:link_weightingI})) for (a), (b), and (c). Link weighting II (Eq.~(\ref{eq:link_weightingII})) for (d), (e), and (f). Hysteresis procedures (a), (d) with $K_i$ from a truncated Gaussian distribution of Eq.~(\ref{eq:tr_gauss_dist}) with $s=1.0$ and $\beta=0.45\pi$; (b),(e) with $K_i$ from a uniform distribution of Eq.~(\ref{eq:uniform_dist}) with $w=0.8$ and $\beta=0.40\pi$; and (c), (f) with $K_i$ from a truncated power-law distribution of Eq.~(\ref{eq:tr_power_law_dist}) with $\gamma=1.3$ and $\beta=0.15\pi$. The parameters $(s,\beta)$, $(w,\beta)$, and $(\gamma,\beta)$ are chosen to be the same as those in Figs.~\ref{fig:01_trg}(d), \ref{fig:01_uni}(d), and \ref{fig:01_trpl}(c), respectively, for direct comparison, and the lines are the same ones theoretically obtained for the model of Eq.~(\ref{eq:model}). The underlying networks are ER networks generated with connection probability $p=0.08$, resulting in a mean degree $\langle k \rangle$ approximately $p(N-1) = 79.92$. 
	For each figure, $R$ values from 5 hysteresis procedures are displayed, each with different $\{K_i\}$, ER network realizations, initial conditions, and, for weighting scheme II, different $\xi$ matrices. For other details, see Fig.~\ref{fig:01_trg} and the text of Subsubsec.~\ref{sbsbsec:ERnet}.
	}
\label{fig:01_wnet}
\end{figure*}

\begin{figure*}[tb!]
    \centering
    \includegraphics[width=17cm]{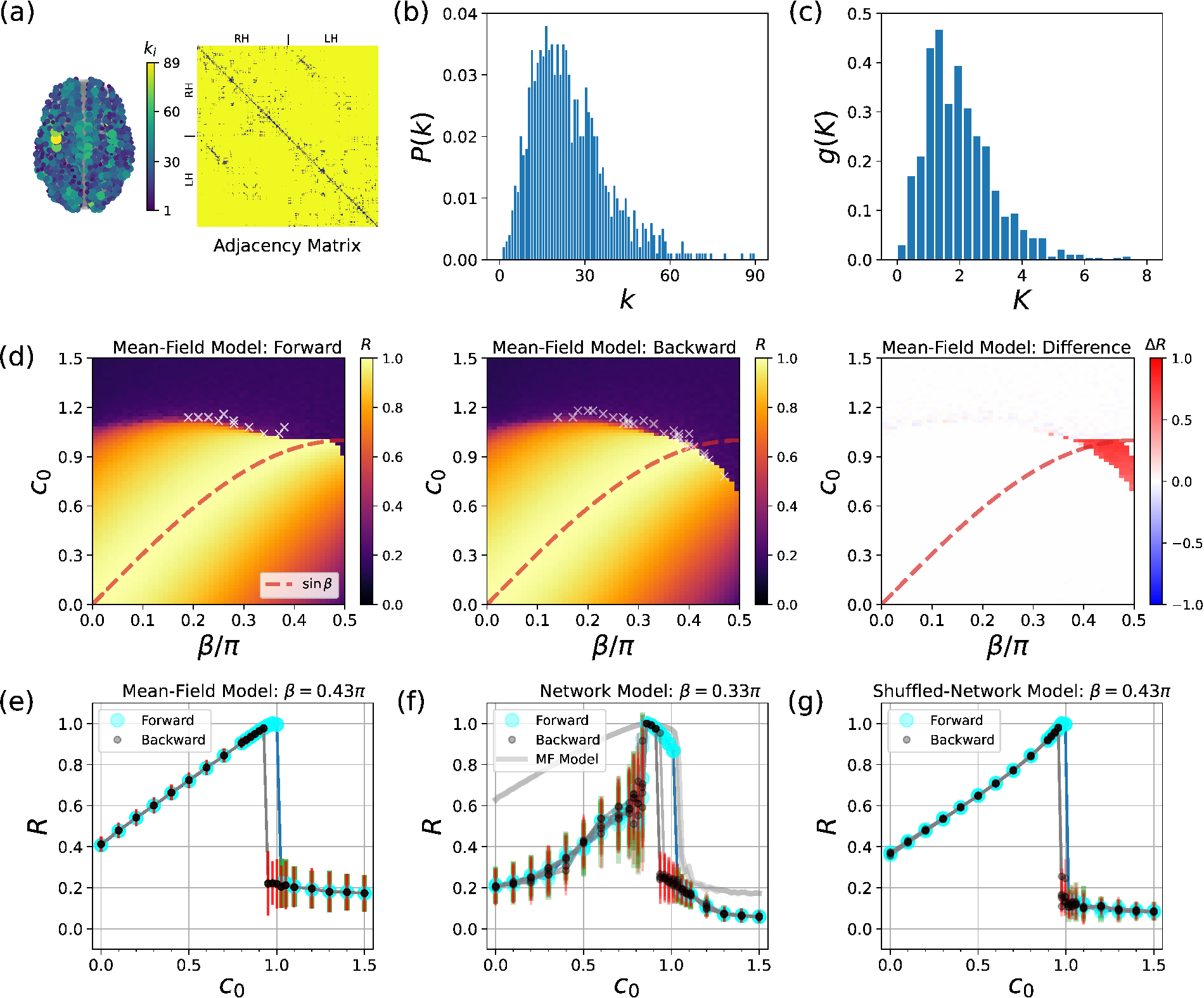}
	\caption{(a) An undirected binary brain network with 1000 nodes ($N=1000$). Left: The network is visualized with nodes located at the positions of the corresponding brain regions. Node sizes are proportional to their degree $k_i$. Edges between nodes are drawn in gray. 
	Right: Adjacency Matrix $A$ for the network.  
	A dot for $A_{ij}=1$ (connected nodes), while no dot for $A_{ij}=0$ (unconnected nodes). 
	RH and LH stand for right hemisphere and left hemisphere, respectively.
	(b) The degree distribution $P(k)$. $P(k) \equiv {n_k}/N$, where $n_k$ is the number of nodes with degree $k$.
	(c) Density histogram $g(K)$ for coupling strength $K$ with bins in the range $[0,9]$ with bin width $0.3$. $K_i \equiv K_m{k_i}/{\langle k \rangle}$, where $\langle k \rangle$ is the average degree and $K_m$ is the mean coupling strength set to $2$ for the comparison with truncated Gaussian and uniform distributions.
	(d) Mean-field model of Eq.~(\ref{eq:model}) with $K_i$ from (c): Hysteresis search in $\beta$-$c_0$ plane using hysteresis procedures for every value of $\beta \in [0,~0.49\pi]$ separated by $\Delta \beta=0.01\pi$. (Forward) $R_\text{forward}(\beta,c_0)$, time-averaged order parameter values for the forward path, averaged over 5 hysteresis procedures with different initial conditions; (Backward) $R_\text{backward}(\beta,c_0)$, $R$ values for the backward path, averaged over the 5 hysteresis procedures; (Difference) $\Delta R(\beta,c_0) \equiv R_\text{forward}(\beta,c_0) - R_\text{backward}(\beta,c_0)$. White X marks indicate data points where the std of the order parameter time series $r(t)$ exceeds $0.15$ for at least one simulation out of $5$. Hysteresis cases (e) in the mean-field model of (d) with $\beta=0.43\pi$, (f) in the network model of Eq.~(\ref{eq:model_network}) with $W_{ij} = K_m A_{ij}/{\langle k \rangle}$ and $\beta=0.33\pi$, where $A$ is the adjacency matrix for the brain network of (a), and (g) in the shuffled-network model of Eq.~(\ref{eq:model_network}) with $W_{ij} = K_m {\tilde A}_{ij}/{\langle k \rangle}$ and $\beta = 0.43\pi$, where $\tilde A$ is the adjacency matrix generated by degree-preserving shuffling of $A$. 
	For each figure of (e), (f) and (g), the lines connecting symbols are intended only as visual guides, and $R$ values from 5 hysteresis procedures are displayed, each with different initial conditions, and, for (g), different shuffled adjacency matrices. Green and red vertical lines are the errorbars denoting the standard deviations of $r$ time series for forward and backward procedures, respectively. In (g), the gray line labeled as `MF Model' represents the $R$ values for the hysteresis procedure from (d) with $\beta=0.33\pi$. For other details, see Figs.~\ref{fig:01_trg}, \ref{fig:02_trg} and the text of Subsubsec.~\ref{sbsbsec:bnet}.
	}
\label{fig:0102_bnet}
\end{figure*}

For the simulations designed to explore the existence of hysteresis, we gradually increase $c_0$ from 0 to 1.5. We refer to this path of increasing $c_0$ values as the {\it forward path}. An incoherent state where the phase of the oscillators are randomly selected from $[0,2\pi)$ is used as the initial condition for $c_0=0$.  For each value of $c_0$, we wait for the transient behavior to subside and measure the time average of $r(t)$, before proceeding to the next $c_0$ value, using the final state from the previous $c_0$ value as the initial condition for the next simulation.  When the maximum value of $c_0=1.5$ is reached, this procedure is mirrored as we gradually decrease $c_0$ from 1.5 to 0, maintaining the approach of using the final state of one simulation as the starting point for the next ({\it backward path}). 
The system typically exhibits a stationary state denoted by the order parameter $r(t)$ with small fluctuation after the transients. However, near the transition to incoherent states and in a small region around $\beta = 0.5\pi$ it shows large fluctuation or large oscillations. 
For the simulations of hysteresis curves ($c_0$-$R$ curves) of Figs. \ref{fig:01_trg}, \ref{fig:01_uni}, \ref{fig:01_trpl}, \ref{fig:01_wnet}, and \ref{fig:0102_bnet}, we use $c_0$ step size $\Delta c_0=0.1$ mainly, but smaller step sizes of $c_0$ are used for ranges where $R$ changes rapidly, near transition points, or inside or near the bistable region. Inside the bistable region, some stable states may not be obtained if the preceding state during the hysteresis procedure falls outside the basin of attraction for those states. Using smaller step size can make the preceding states more similar to these stable states, increasing the likelihood that the preceding state lies within the basin of attraction.  
There still can be earlier transitions to the other stable states near transition points obtained theoretically. This is due to the finite size effect and the sensitivity to initial conditions for the states in those regions.
In simulations for the systematic searches for the hysteresis in $\beta$-$c_0$ plane as in Fig. \ref{fig:02_trg}, $c_0$ step size is fixed to $\Delta c_0=0.02$ to provide a fine resolution. 

%===================================
% VVV Truncated Gaussian Dist VVVVVV 
%===================================

\subsection{\label{subsec:tr_gaussian_dist}Truncated Gaussian coupling strength distributions}
First, we investigate the system with a truncated Gaussian distribution given by
\begin{eqnarray}
	g(K)&=&C_g {\rm e}^{-\frac{1}{2}(K-K_m)^2/s^2}~~~\text{for $K \in [K_a, K_b]$},
\label{eq:tr_gauss_dist}
\end{eqnarray}
where $C_g$ is the normalization factor, $K_a \in (0, K_m]$, and $K_b = 2K_m-K_a$. 
The lower truncation is introduced to ensure all sampled $K$ values from the distribution are positive, which is necessary for our coupling strengths.
In this symmetric setting for the distribution, $K_m$ is the average value of $K$. 
With the fixed values of $K_m$, $K_a$, and $K_b$, the parameter $s$ controls the width of the distribution (Fig. \ref{fig:01_trg}(a)). 
This parameterization allows us to systematically change the shape of the distribution and explore the effects of it on the hysteresis of the system.

For the simulations and analysis, we use $K_m=2.0$, $K_a=0.01$ and $K_b=3.99$. We use ten different sets of initial condition and $\{K_i\}$ for the simulations of Figs.~\ref{fig:01_trg}(c) and (d).

In Fig.~\ref{fig:01_trg}(a), we illustrate the truncated Gaussian coupling strength distributions for $s=0.5$, $s=1.0$, and $s=3.0$, respectively. As $s$ increases, the peak of the distribution broadens and flattens.  In Fig.~\ref{fig:01_trg}(b), we present the critical threshold $c_0^*$ for uniformly incoherent states as a function of $\beta$ for the three different values of $s$. For a given value of $s$ and $\beta$, the system can exhibit a uniformly incoherent state when $c_0 > c_0^*$. The curve for each value of $s$ shows monotonically decreasing behavior with increasing $\beta$ for this coupling strength distribution.  
With increasing $s$, the critical threshold $c_0^*$ decreases for the same value of $\beta$. This indicates that the range of $c_0$ for which a uniformly incoherent state is obtained extends further into smaller values of $c_0$. The curve for $c_0 = \sin\beta$, with which the system has  stable fully locked states with $r=1$, thus $R=1$, is shown for comparison. This curve serves as a reference relative to which we can consider things such as $c_0^*$ curves and hysteresis ranges for $c_0$.

Figures~\ref{fig:01_trg}(c) and (d) display cases without and with hysteresis, respectively. Time-averaged order parameter values ($R$) obtained from numerical simulations, denoted by symbols, match well with the order parameter $r$ values of at least linearly neutrally stable states predicted by the self-consistency analysis and stability analysis. Solid red lines represent at least neutrally stable states, while dashed gray lines indicate unstable states. 
Note that each symbol corresponds to a $R$ value from a single simulation at a specific $c_0$ value. 
With $s=0.5$ and $\beta=0.45\pi$, as $c_0$ increases from $0$ to $\sin\beta$ through the forward path, there is a gradual increase in $R$, from approximately 0.6 to 1 (Fig. \ref{fig:01_trg}(c)). When $c_0$ surpasses $\sin\beta$, $R$ begins to decrease, eventually stabilizing at a low level. Beyond this point, further increases in $c_0$ do not significantly affect $R$. The backward path of decreasing $c_0$ from $1.5$ to $0$ gives similar values of $R$. There is no hysteresis for this case.
In contrast, the system shows abrupt transitions and hysteresis with $\beta=0.45\pi$ and wider coupling strength distribution given by $s=1.0$ (Fig. \ref{fig:01_trg}(d)). In this case, there are sudden jumps from a coherent state to an incoherent state, and vice versa, as $c_0$ changes.
The system exhibits two distinct states, one with high $R$ and the other with low $R$, for the same values of $c_0$ in the forward and backward paths, respectively. Figures~\ref{fig:01_trg} (e) and (f) show the phases of individual oscillators of the two states with $c_0$ in the bistable range. Figure~\ref{fig:01_trg}(e) depicts the partially locked state with high $R$, characterized by drifting oscillators with low $K_i$ values and locked oscillators with high $K_i$ values, while Fig.~\ref{fig:01_trg}(f) shows the uniformly incoherent state with low $R$. 

In Fig.~\ref{fig:02_trg}, for each $\beta$ value in $[0, 0.49\pi]$ separated by an interval of $\Delta \beta = 0.01\pi$, we examine the presence of hysteresis by performing a hysteresis procedure. 
For a given $s$ value, we plot $R$ values for the forward path in the panel labeled Forward, and $R$ values for the backward path in the panel labeled Backward, respectively. The difference $\Delta R$ values, defined as $\Delta R(\beta,c_0) \equiv R_\text{forward}(\beta,c_0) - R_\text{backward}(\beta,c_0)$, are plotted in the panel labeled Difference. Here, $R_\text{forward}(\beta,c_0)$ and $R_\text{backward}(\beta,c_0)$ are $R$ values obtained along the forward and backward paths for $c_0$ with a fixed value of $\beta$, respectively. Non-zero $\Delta R$ values indicate the presence of hysteresis, with larger absolute values suggesting stronger hysteresis effects.

In Fig.~\ref{fig:02_trg}(a), for the smallest value of $s$, $s=0.5$, corresponding to the least inhomogeneity in $K_i$, 
the system exhibits states with high $R$ for most of $(\beta,c_0)$ values except those in the upper-right corner of the $\beta$-$c_0$ plane (high $\beta$, high $c_0$) for both of the forward and backward paths. We observe significant hysteresis, evidenced by high $\Delta R$ in red, in the narrow region bounded by $\beta=0.49\pi$, bounded below by $c_0 = c_0^*$ (cyan solid line, representing the critical threshold for uniformly incoherent states), and bounded above by $c_0 = \sin\beta$ (red dashed line). In this region, the system exhibits bistability between a uniformly incoherent state and a partially (fully) locked state. 
Because of the early transitions mentioned in the introductory part of Sec.~\ref{sec:Simulation}, the lower boundary of the red region does not match perfectly with the $c_0^*$ curve. 
To examine abrupt phase transitions through the steepness of the $c_0$-$R$ curve in the range from the peak with $r=1$ to the transition point, for the forward procedure, we numerically search for two key boundaries in the $c_0$ range $[\sin\beta,1.5]$ using the self-consistency equation (Eq.~(\ref{eq:self_consistency})) : $c_0^{0.9}$ and $c_0^{0.1}$. Here, $c_0^{0.9}$ is defined as the maximum $c_0$ for states with $r\geq 0.9$, provided that $r<0.9$ for $c_0=1.5$. If $r\geq0.9$ for $c_0=1.5$, then $c_0^{0.9}$ does not exist for that particular case. $c_0^{0.1}$ represents the minimum $c_0$ for near incoherent states with $r\leq 0.1$. The proximity of these two boundaries indicates a steeper $c_0$-$R$ curve. We observe an abrupt transition to an incoherent state in this region.

As we increase the value of $s$ to $s=1.0$ and to $s=3.0$ in Fig.~\ref{fig:02_trg}(b) and (c), respectively, 
the region of low $R$ values expand near to the $c_0^*$ curves for the forward paths and the region of hysteresis expands to include lower values of $\beta$ and $c_0$. The lower boundaries of the regions match well with the curves of $c_0^*$ whose value decreases for the same value of $\beta$ with increasing $s$. The hysteresis here involves only the bistability between a uniformly incoherent state and a partially (fully) locked state as in Fig.~\ref{fig:02_trg}(a). Note that the hysteresis regions extend to the $c_0$ values above $c_0=\sin\beta$ for lower $\beta$ values, in contrast to the case of Fig.~\ref{fig:02_trg}(a).  

This subsection shows that coupled identical oscillators with a more realistic coupling function $H(\theta) = c_0 + \sin(\theta-\beta)$ and coupling strength inhomogeneity can exhibit abrupt phase transitions and hysteresis. With truncated Gaussian distributions for coupling strengths, the hysteresis regions lie near $\beta=\pi/2$, and are approximately bounded above by $c_0=\sin\beta$ curve and bounded below by the $c_0^*$ curves, above which uniformly incoherent states are observed. The hysteresis involves only the bistability between a partially (fully) locked state and a uniformly incoherent state.  As the distribution broadens, the hysteresis region expands to include lower values of $\beta$ and $c_0$. 
%===================================
% ^^^ Truncated Gaussian Dist ^^^^^^
%===================================

%===================================
% VVV ---- Uniform Dist ----- VVVVVV
%===================================
%=============================
%Fig. 3
%=============================

%=============================
%Fig. 4
\begin{figure*}[tb!]
    \centering
    \includegraphics[width=17cm]{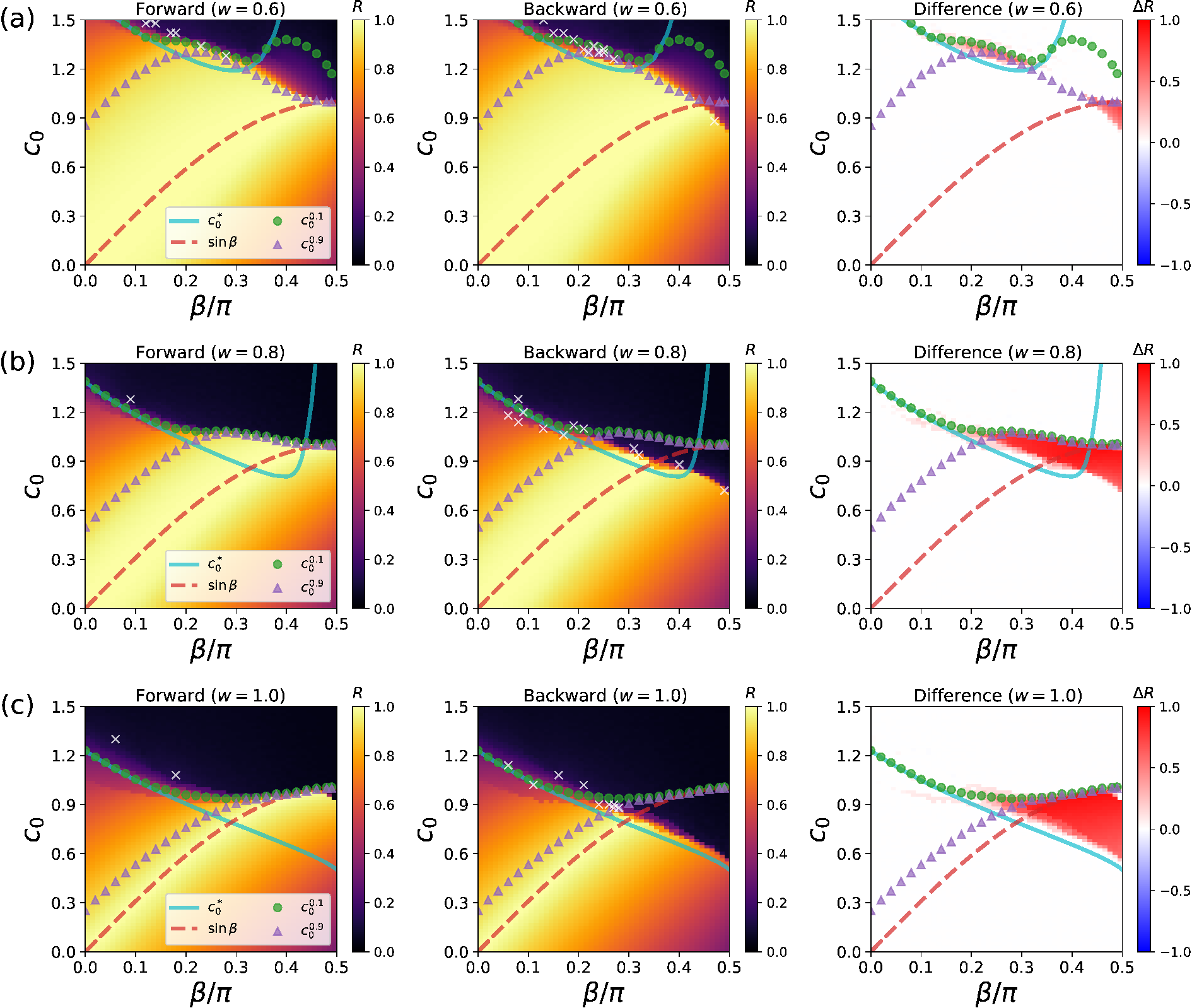}
	\caption{Hysteresis search in $\beta$-$c_0$ plane using hysteresis procedures with uniform distributions for every value of $\beta \in [0,~0.49\pi]$ separated by $\Delta \beta=0.01\pi$.  (a) $w=0.6$: (Forward) $R_\text{forward}(\beta,c_0)$, time-averaged order parameter values for the forward path, averaged over 5 hysteresis procedures with different $\{K_i\}$ and initial conditions; (Backward) $R_\text{backward}(\beta,c_0)$, $R$ values for the backward path, averaged over the 5 hysteresis procedures; (Difference) $\Delta R(\beta,c_0) \equiv R_\text{forward}(\beta,c_0) - R_\text{backward}(\beta,c_0)$. (b) $w=0.8$. (c) $w=1.0$. Other details are as in Fig.~\ref{fig:02_trg}.  
	}
    \label{fig:02_uni}
\end{figure*}
%=============================

%=============================
%Fig. 5
%=============================

%=============================
%Fig. 6
\begin{figure*}[tb!]
    \centering
    \includegraphics[width=17cm]{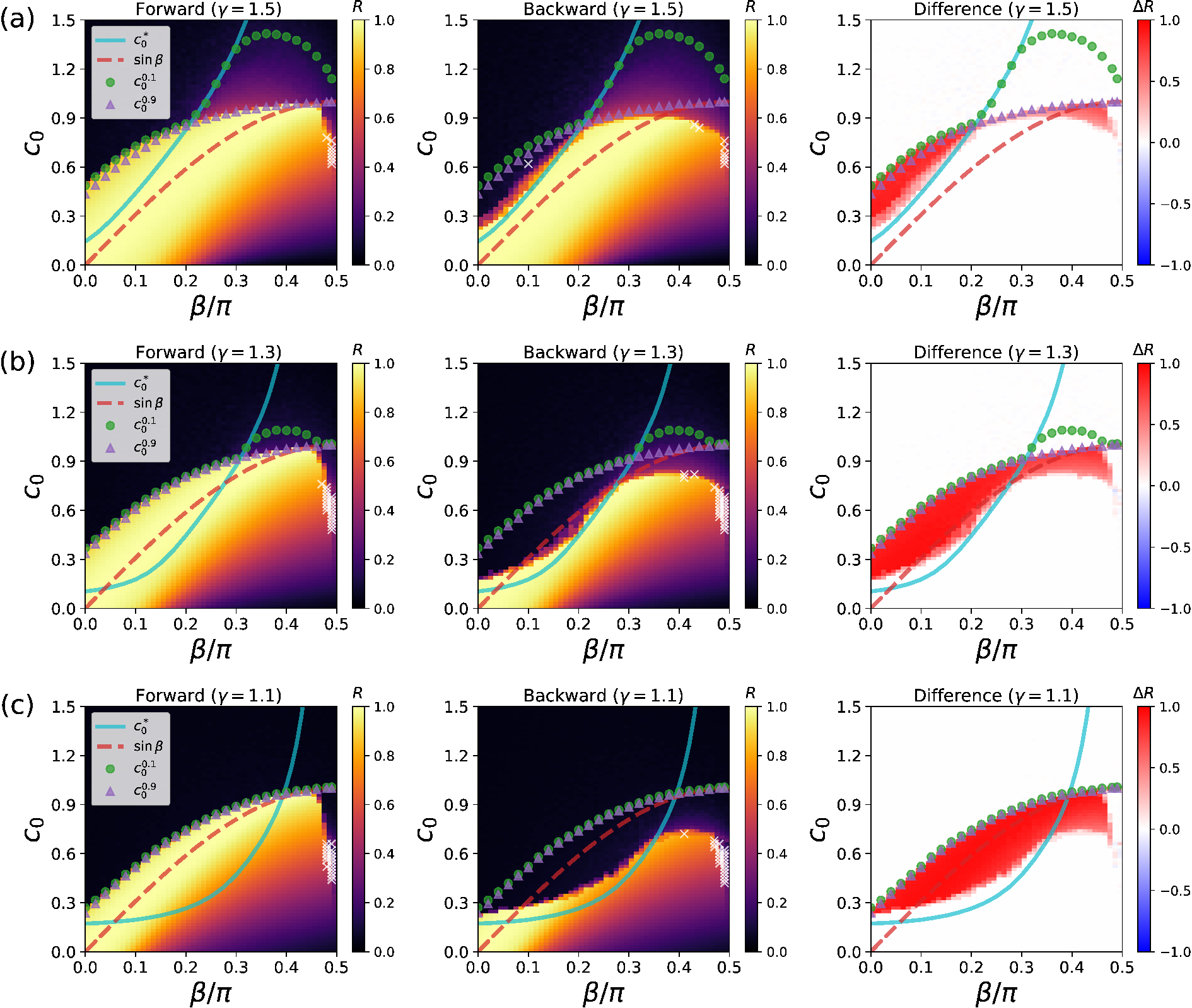}
	\caption{Hysteresis search in $\beta$-$c_0$ plane using hysteresis procedures with truncated power-law distributions for every value of $\beta \in [0,~0.49\pi]$ separated by $\Delta \beta=0.01\pi$.  (a) $\gamma=1.5$: (Forward) $R_\text{forward}(\beta,c_0)$, time-averaged order parameter values for the forward path, averaged over 5 hysteresis procedures with different $\{K_i\}$ and initial conditions; (Backward) $R_\text{backward}(\beta,c_0)$, $R$ values for the backward path, averaged over the 5 hysteresis procedures; (Difference) $\Delta R(\beta,c_0) \equiv R_\text{forward}(\beta,c_0) - R_\text{backward}(\beta,c_0)$. (b) $\gamma=1.3$. (c) $\gamma=1.1$. Other details are as in Fig.~\ref{fig:02_trg}.  
	}
\label{fig:02_tr_power_law}
\end{figure*}
%=============================

%=============================
%Fig. 7
%=============================

%=============================
%Fig. 8
%=============================

\subsection{\label{subsec:uniform_dist}Uniform coupling strength distributions}
Next, we investigate the system with a uniform distribution given by
\begin{eqnarray}
g(K)&=& C_u~\text{for $K \in [K_m(1- w), K_m(1+w)]$,} \nonumber \\
        &&~~~~~\text{with~} w \in [0,1],
\label{eq:uniform_dist}
\end{eqnarray}
where $C_u$ is the normalization factor, and $K_m$ is the average value of $K$. With a fixed value of $K_m$, the parameter $w$ controls the width of the distribution: the larger $w$, the larger the width of the distribution. 
For the simulations and analysis, we use $K_m=2.0$.

In Figs.~\ref{fig:01_uni}(a) and (b), we plot the coupling strength distributions and  the critical threshold $c_0^*$ for uniformly incoherent states as a function of $\beta$ for $w=0.6$, $w=0.8$, and $w=1.0$, respectively. 
The $c_0^*$ curve for $w=1.0$ shows monotonically decreasing behavior with increasing $\beta$ as with the truncated Gaussian coupling strength distribution. In contrast, the curves for $w=0.6$ and $w=0.8$ show non-monotonic behaviors. They monotonically decrease first and then monotonically increase. 
With increasing $w$, the critical threshold $c_0^*$ decreases for the same value of $\beta$, indicating that the range of $c_0$ for which a uniformly incoherent state is linearly neutrally stable increases with larger $w$.

Figures~\ref{fig:01_uni}(c)-(f) display cases with hysteresis. As in the cases with truncated Gaussian distributions in Figs.~\ref{fig:01_trg} and \ref{fig:02_trg}, the same kind of bistability between a partially (fully) locked state and a uniformly incoherent state is observed in Figs.~\ref{fig:01_uni}(c)-(e). The lower bounds of the bistability regions are given by $c_0^*$. Meanwhile, a different type of bistability between a partially (fully) locked state with high $R$ and a partially locked state with low $R$ are observed for $c_0 < c_0^*\approx 2.09$ range where uniformly incoherent states are unstable in Fig.~\ref{fig:01_uni}(f). 
In this figure, the red line for $c_0 \in [1.0,1.5]$ represents partially locked states with very small $R$, rather than uniformly incoherent states. 
Figures~\ref{fig:01_uni}(g) and (h) show the partially locked states with high and low $R$, respectively, which constitute the bistable pair observed along the forward and backward paths of (f).
Note that the partially locked state of Fig.~\ref{fig:01_uni}(h) is with oscillators with small $K_i$ locked contrary to the partially locked states of Figs.~\ref{fig:01_trg}(e) and \ref{fig:01_uni}(g) where oscillators with large $K_i$ are locked. 

Figure~\ref{fig:02_uni} shows hysteresis search in $\beta$-$c_0$ plane for the cases with uniform coupling strength distributions. In Fig.~\ref{fig:02_uni}(a), for the smallest value of $w$, $w=0.6$, we observe significant hysteresis in two separate regions: the narrow region above ${c_0}^*$ curve (cyan solid line) and the triangular region bounded by $c_0=\sin\beta$ (red dashed line) and $\beta=0.49\pi$. The former region is far from the $c_0=\sin\beta$ curve and covers high $c_0$ and low $\beta$ values not covered by the cases with truncated Gaussian distributions in Fig.~\ref{fig:02_trg}. In this region, bistability between a partially locked state and a uniformly incoherent state is observed along the hysteresis procedures.  
In the latter region, hysteresis involving bistability between a partially (fully) locked state with high $R$ and a partially locked state with low $R$ as observed in Fig.~\ref{fig:01_uni}(f) occurs. 

As we increase the value of $w$ to $w=0.8$ in Fig.~\ref{fig:02_uni}(b), the two hysteresis regions merge into one. This hysteresis region is composed of the area above the $c_0^*$ curve and the area below it. The upper bound of the hysteresis region is approximately given by $c_0^{0.9}$ boundary, which lies very close to $c_0^{0.1}$ boundary in this case. 
In the hysteresis region above the $c_0^*$ curve, the system exhibits bistability between a partially (fully) locked state and a uniformly incoherent state.
In the hysteresis region below the curve, there is bistability between a partially (fully) locked state with oscillators having large $K_i$ locked and a partially locked state with oscillators having small $K_i$ locked.  With $w=1.0$ in Fig.~\ref{fig:02_uni}(c), the $c_0^*$ curve is monotonically decreasing, and the hysteresis region is entirely above this curve. The distribution with $w=1.0$ resembles the truncated Gaussian distribution with $s=3.0$, and the hysteresis regions for both cases look very similar. In this case, the hysteresis characterized only by the bistability between a uniformly incoherent state and a partially (fully) locked state.  

This subsection shows that there are cases where $c_0^*$ curves exhibit non-monotonic behaviors with uniform distributions for coupling strengths, in contrast to the monotonic decreasing behaviors with truncated Gaussian distributions. The curves monotonically decrease first and then monotonically increase as $\beta$ increases. In these cases, we observe hysteresis in regions lying both above and below the $c_0^*$ curves, while the hysteresis regions lie only above the curves with truncated Gaussian distributions. Notably, the hysteresis observed below the $c_0^*$ curves involves the bistability between a partially (fully) locked state with oscillators having large $K_i$ and a partially locked state with oscillators having small $K_i$ locked.
%===================================
% ^^^ ---- Uniform Dist ----- ^^^^^^
%===================================

%===================================
%---Truncated Power-law Dist ----VVV 
%===================================
\subsection{\label{subsec:tr_power_law_dist}Truncated Power-law coupling strength distributions}
Finally, we investigate the system with a truncated power-law distribution given by
\begin{eqnarray}
	g(K)&=& C_s K ^{-\gamma}~~~\text{for $K \in [K_a, K_b]$,} \nonumber \\
        &&~~~~~~~~~~~~~\text{with~} \gamma >0,
\label{eq:tr_power_law_dist}
\end{eqnarray}
where $C_s$ is the normalization factor and $K_b>K_a > 0$.
The upper truncation is introduced to ensure numerical stability of the simulations. 
With the fixed values of $K_a$ and $K_b$, the parameter $\gamma$ controls the shape of the distribution. Larger $\gamma$ values result in a steeper decline in the distribution.
For the simulations and analysis, we use $K_a=0.01$ and $K_b=10$.

In Fig.~\ref{fig:01_trpl}(a), we show the coupling strength distributions for $\gamma=1.5$, $\gamma=1.3$, and $\gamma=1.1$, respectively. 
In Fig.~\ref{fig:01_trpl}(b), we display the critical threshold $c_0^*$ for uniformly incoherent states as a function of $\beta$ for the three different values of $\gamma$. 
In contrast to the previous two types of coupling strength distributions, the $c_0^*$ curves monotonically increase with $\beta$ for the truncated power-law distributions. The $c_0^*$ values are very small around $c_0 = \sin\beta$ for small $\beta$ values and then increase with $\beta$. The $c_0^*$ values decrease with $\gamma$ decrease for a fixed value of $\beta$ except in a small range near $\beta=0$. This means the range for uniformly incoherent states extends to lower $c_0$ values for most values of $\beta$ as we decrease $\gamma$. 
Given the characteristics of the $c_0^*$ curves, the system exhibits hysteresis for smaller values of $\beta$ and $c_0$ compared to the cases with the previous two distributions (Figs.~\ref{fig:01_trpl}(c) and (d)). 
For the hysteresis observed with $c_0 > c_0^*$ as in Figs.~\ref{fig:01_trpl}(c) and (d), the system exhibits bistability between a partially locked state and a uniformly incoherent state (shown in Figs.~\ref{fig:01_trpl}(e) and (f), (g) and (h), respectively).
Furthermore, similar to the cases with uniform distributions, for the hysteresis below the $c_0^*$ curves, the system exhibits bistability between two partially locked states: one where oscillators with large $K_i$ are locked, and another where oscillators with small $K_i$ are locked.

Figure~\ref{fig:02_tr_power_law} shows the results of hysteresis search in $\beta$-$c_0$ plane using hysteresis procedures with truncated power-law distributions. Compared to the cases with truncated Gaussian and uniform distributions, the region of low $R$ values occupies a significantly larger portion of the plane, including areas with low $\beta$ and $c_0$ in the panels labeled Forward. As with uniform distributions, we observe significant hysteresis both in the regions above and below the $c_0^*$ curves. The hysteresis occurs for most of the $\beta$ values, including $\beta=0$, and the $c_0$ value with which hysteresis is observed has a minimum less than $0.3$ for $\beta=0$. The hysteresis regions exhibit distinctly different shapes compared to those with truncated Gaussian and uniform distributions.

Based on the transformed model of Eq.~(\ref{eq:model_tr}), our observation of hysteresis in the positive $c_0$ range for $\beta=0$ is consistent with the previous studies of systems of coupled nonidentical oscillators with the coupling function $H(\theta) = \sin\theta$~\cite{gomez2011,xzhang2013}. These studies observed abrupt phase transitions and hysteresis in the presence of positive degree-frequency correlation~\cite{gomez2011} or coupling strength-frequency correlation~\cite{xzhang2013} when degree or coupling strength follows a power-law distribution. 
%---

As mentioned in Sec.~\ref{sec:Model}, the transformed model has the same form of the Sakaguchi-Kuramoto model with the coupling function $H(\theta) = \sin(\theta-\beta)$ for coupled nonidentical oscillators on networks with a positive degree-frequency correlation \cite{kundu2017}.
While we observe hysteresis for most $\beta$ values, Ref. \cite{kundu2017} reported hysteresis with uncorrelated scale-free networks showing power-law degree distributions (exponents $\gamma = 2.7$ and $2.8$) only for small $\beta$ values below a threshold. This discrepancy might stem from differences in distributions and $\beta$ ranges studied. Our coupling strength distributions use exponents $\gamma = 1.1$, $1.3$, and $1.5$, while they limited their investigation to $\beta \leq 1.2 (\approx 0.38\pi)$. We observe that as $\gamma$ increases, $c_0^*$ tends to move upward and leftward. Consequently, the hysteresis region above the $c_0^*$ curve shrinks, while the region below converges toward $\beta = \pi/2$ (Fig.~\ref{fig:02_tr_power_law}(c) $\rightarrow$(a)). It is possible that Ref.~\cite{kundu2017} observed only the hysteresis region corresponding to that above the $c_0^*$ curve in our model.

%---
This subsection demonstrates that hysteresis can also occur in systems with truncated power-law coupling strength distributions. In this case, we observe hysteresis across a wide range of $\beta$ values, from $0$ to those close to $\pi/2$. This contrasts significantly with the hysteresis regions observed in systems with truncated Gaussian or uniform distributions, as discussed in previous subsections. The extensive range of hysteresis in the truncated power-law case highlights a distinct characteristic of these systems, suggesting that the nature of the coupling strength distribution plays a crucial role in determining the hysteresis behavior.

\subsection{Models with weighted networks for connectivity between oscillators}
In this subsection, to confirm that the results about the hysteresis of the previous subsections are valid in more realistic situations where oscillators are connected according to a weighted network, we study the following network model of Eq. (\ref{eq:model}). 
\begin{eqnarray}
    \frac{d \theta_i}{dt} &=& \omega + \sum_{j=1}^N W_{ij}\big [c_0+\sin(\theta_j - \theta_i-\beta)\big],  
\label{eq:model_network} \\
	&&i=1,2,...,N, ~c_0 \geq 0, ~\beta \in [0,\pi/2), ~W_{ij}\geq 0, \nonumber
\end{eqnarray}
where $W_{ij}$ is the element of the weighted adjacency matrix $W$ giving the link strength when oscillator $j$ influencing oscillator $i$. 
Note that in general $W_{ij}$ is asymmetric meaning the strength of $j$ influencing $i$ is different from that of $i$ influencing $j$. 

Through the same transformation used for Eq. (\ref{eq:model_tr}), Eq. (\ref{eq:model_network}) can be written as
\begin{eqnarray}
    \frac{d \phi_i}{dt} &=& K_i+\frac{1}{c_0}\sum_{j=1}^N W_{ij}\sin(\theta_j - \theta_i-\beta),  
\label{eq:model_network2} 
\end{eqnarray}
where $K_i=\sum_{j=1}^N W_{ij}$. We are interested in the cases with the values of $K_i$ is inhomogeneous. As in the case of the original model of Eq. (\ref{eq:model}), this transformed model has the form of coupled nonidentical oscillators with a positive correlation between intrinsic frequencies and coupling strengths \cite{xzhang2013}.

\subsubsection{\label{sbsbsec:ERnet}Weighted Erd{\H o}s-R{\'e}nyi random network for connectivity}
To compare the results from the models of Eq.~(\ref{eq:model_network}) with those of previous subsections, we consider the cases with $\sum_{j=1}^N W_{ij} = K_i$, where the coupling strength distributions of the subsections \ref{subsec:tr_gaussian_dist}, \ref{subsec:uniform_dist}, and \ref{subsec:tr_power_law_dist} are used for randomly generating $\{K_i\}$.
For the connectivity representing the presence of connections between oscillators, we use Erd{\H o}s-R{\'e}nyi (ER) random networks, which are simply generated by connecting nodes with a probability \cite{erdos1959,barabasi2016}. 
Since intrinsic structural properties such as clustering and modular structures observed in real-world networks \cite{barabasi2016} can affect the synchronization \cite{arenas2008} and hysteresis behaviors, we choose ER random networks, which do not have such properties, for the generalization of the mean-field model of Eq.~(\ref{eq:model}).      
We generate ER random networks with connection probability $p$, and then make them weighted by applying one of two weighting schemes to the links: I) Equally distributing $K_i$ to the $k_i$ in-coming links of oscillator $i$:
\begin{eqnarray}
W_{ij} = \left (\frac{K_i}{k_i} \right) A_{ij},
\label{eq:link_weightingI}
\end{eqnarray}
where $A_{ij}$ is the adjacency matrix of the ER network, where $A_{ij}=1$ if $i$ and $j$ are connected and $A_{ij}=0$ otherwise. Here, $k_i$ is the in-degree of node $i$ meaning there are $k_i$ in-links to $i$, and $K_i$ is the coupling strength as in the previous sections.
II) Unequally distributing $K_i$ to the $k_i$ in-coming links of oscillator $i$ according to:
\begin{eqnarray}
W_{ij} = \left (\frac{\xi_{ij} K_i}{\sum_{j\in \mathcal{N}_i} \xi_{ij}} \right) A_{ij},
\label{eq:link_weightingII}
\end{eqnarray}
where $\xi_{ij}$ are randomly assigned from $(0,1]$ for weighting, introducing heterogeneity in link strengths, and $\mathcal{N}_i$ is the set of oscillators influencing $i$.
For both of the weighting schemes, the link strengths $W_{ij}$ are set to satisfy $\sum_{j=1}^N W_{ij} = K_i$ for all $i$.

Figure~\ref{fig:01_wnet} displays the hysteresis curves obtained with the network models using the two link weighting schemes when the connection probability $p$ for the generation of ER random networks is $0.08$. For the three distributions for $K_i$, the curves closely match the theoretically obtained values for the mean-field model of Eq.~(\ref{eq:model}), denoted by the lines. Larger values of $p$ do not change the results significantly. In contrast, with smaller values of $p$, the $R$ values are lowered and the $c_0$ ranges for the hysteresis are narrowed. 

These results show that the generalized model with a weighted random network for the connectivity between oscillators can exhibit abrupt phase transitions and hysteresis matching those of the mean-field model when there is a sufficiently large number of connections. This confirms the finding of the preceding subsections that the simplified realistic coupling function and coupling strength inhomogeneity can cause the abrupt phase transitions and hysteresis, and generalize it to the cases of systems with weighted connectivities.   

\subsubsection{\label{sbsbsec:bnet}A brain network for connectivity}
Finally for the studies with networks, we investigate hysteresis in the network model of Eq.~(\ref{eq:model_network}) using a human brain network for the connectivity as an example for the cases with real-world complex networks. 
Specifically, we use a structural connectivity map of human cerebral cortex comprising $N=1000$ cortical regions constructed from diffusion spectrum imaging (DSI) \cite{shafiei2019,shafiei2019data}. In this human brain network, each node corresponds to a brain region, and an edge between two nodes represents the existence of neural connections between the corresponding regions (Fig. \ref{fig:0102_bnet}(a)). The network is undirected due to the tracing method for neural connections based on DSI. Here, focusing only on the connectivity, we consider an unweighted network. 
The degrees of the brain regions range between $1$ and $89$, as shown in Fig. \ref{fig:0102_bnet}(b), reflecting the heterogeneity of connectivity in the human brain. This inhomogeneity in the degrees has the role of the coupling strength inhomogeneity.

Let us first consider the mean-field model of Eq.~(\ref{eq:model}) with $K_i = K_m k_i/{\langle k \rangle}$ where $k_i$ is the degree of node $i$ of the brain network, ${\langle k \rangle}$ is the mean degree, and $K_m$ is the mean coupling strength set to $2$ for the comparison with truncated Gaussian and uniform distributions. Other values of $K_m >0$ do not change the results qualitatively.   
The density histogram $g(K)$ shown in Fig.~\ref{fig:0102_bnet}(c) resembles more like a truncated Gaussian distribution but exhibits positive skewness, with the distribution tail extending further to the right of the peak compared to the left side.  
Figure~\ref{fig:0102_bnet}(d) displays the results from hysteresis search in $\beta$-$c_0$ plane. The hysteresis occurs in the region near $\beta=\pi/2$ and $c_0=\sin\beta$ curve similarly as in the case with truncated Gaussian distribution with $s=1.0$ (Fig.~\ref{fig:02_trg}(b)). 
In contrast to the case with truncated Gaussian distribution with $s=1.0$, the region with high $R$ values is reduced significantly for $\beta \lesssim 0.3$ for both of the forward and the backward paths. The hysteresis curve in Fig.~\ref{fig:0102_bnet}(e) has the typical shape observed with other distributions of Figs.~\ref{fig:01_trg}, \ref{fig:01_uni}, and \ref{fig:01_trpl}. These results show that coupling strength distribution derived from real-world systems can induce hysteresis.

Now, we study the network model of Eq.~(\ref{eq:model_network}) with $W_{ij} = K_m A_{ij}/{\langle k \rangle}$, where $A$ is the adjacency matrix for the brain network. As shown in Fig.~\ref{fig:0102_bnet}(f), the system exhibits hysteresis but the hysteresis curve is different from that of the mean-field model. Moreover, the hysteresis is observed in a different region of the $\beta$-$c_0$ plane. The system shows large fluctuations or oscillation-like behavior for each $c_0$ value along the hysteresis procedure, in contrast to the relatively small fluctuations in the mean-field model. 
We think the difference in the shape of the hysteresis curve and the large fluctuations are due to the intrinsic structural properties of the brain network, such as the tendency of the nodes connected to nearby nodes, the presence of modular structures, where connections inside a module are denser than those toward outside from the module, and the connectivity patterns among nodes of varying degrees. 
Each of these intrinsic structural properties can independently contribute to the deviation from the mean-field approximation by introducing local correlations and heterogeneities not captured in the average behavior.
As shown in the adjacency matrix (Fig.~\ref{fig:0102_bnet}(a)) for the brain network, the nodes in each hemisphere of the brain are connected predominantly with those within the same hemisphere, and form a module. We can see more modular structures inside the hemispheres. 
In a study on chimera states, which show spatial coexistence of phase-locked and drifting oscillators, order parameter oscillation was observed in networks with modular structures using coupling functions of the form $H(\theta)=\sin(\theta-\beta)$ \cite{abrams2008}. This is relevant to our observations of large fluctuations in the brain network model. Furthermore, researches on explosive synchronization in coupled nonidentical oscillators have demonstrated that frequency-frequency correlation and degree-degree correlation among connected nodes affect the emergence of explosive transitions and hysteresis \cite{lzhu2013,chen2013,li2013,nadal2015,peron2015}. In our network model, oscillators with high degrees act like ones with large intrinsic frequencies in systems of non-identical oscillators due to $c_0$ in the coupling function and thus frequency-frequency correlation is directly related to degree-degree correlation. We need to consider how the connectivity patterns among nodes of varying degrees affect the abrupt phase transitions and hysteresis in our model in the future works.   

To test whether intrinsic structural properties in the brain network significantly affect the synchronization and the hysteresis behavior, we shuffle the adjacency matrix while retaining the degree sequence, the undirectedness, and the connectedness, and avoiding self-connections \cite{bounova2015}. We then use this shuffled matrix for simulations. Although the $R$ values are slightly lowered and the range of $c_0$ for hysteresis is slightly narrowed, the network model with the shuffled network shows results similar to the mean-field model (Fig.~\ref{fig:0102_bnet}(f)). The slight deviation from the mean-field model is probably due to the nodes with small degrees as we observe such deviation in Subsubsec.~\ref{sbsbsec:ERnet}. The close match between the results from the mean-field model and the network model with the shuffled brain network, which lacks the intrinsic structural properties of the network, validates the idea that the properties affect the synchronization and the hysteresis behavior significantly.

The results of this subsubsection generalize the findings of the preceding subsections by demonstrating that network models of coupled identical oscillators, which use simplified realistic coupling functions and realistic networks with inhomogeneous degree distributions, can exhibit abrupt phase transitions and hysteresis.

\section{Summary and Conclusions}
In summary, we have investigated the existence of hysteresis in a generalized Kuramoto model with identical oscillators coupled with a more realistic coupling function and inhomogeneous coupling strengths which determine how strongly oscillators are coupled to others. 
By retaining the usually ignored constant term in the coupling function, we have shown that the system can be regarded as coupled nonidentical oscillators with inhomogeneous effective intrinsic frequencies which arise from the interplay of the constant term and the coupling strengths.   
Without any special factors explicitly or implicitly incorporated to the model, a positive correlation between the effective intrinsic frequencies and the coupling strengths emerges naturally in the system. This correlation serves as a condition that can lead to abrupt phase transitions and hysteresis.  
Through numerical simulations and analysis using various coupling strength distributions, we have found that this system can exhibit abrupt phase transitions and hysteresis. The distribution of coupling strengths substantially affects the hysteresis regions within the parameter space of the coupling function.
Numerical simulations with a network version of the model, specifically one using a brain network, generalize our findings to the cases of complex real-world systems.

This finding is particularly significant as it demonstrates how a simplified, realistic coupling function and coupling strength inhomogeneity can inherently produce conditions favorable for abrupt phase transitions and hysteresis. This contrasts with previous models that often required explicit additional factors to achieve similar results. Moreover, our study shows that the properties of oscillators and their interactions, which determine the coupling function, can lead to such phenomena when coupling strengths are inhomogeneous.

To further expand on these insights, subsequent work could focus on the studies involving nonidentical oscillators and real complex networks with various properties in addition to the inhomogeneity in the coupling strengths. 
Our results, along with these future directions, can contribute to the understanding of hysteresis and related behaviors in various complex systems, such as the brain \cite{kleinschmidt2002,hkim2018,sepulveda2018}.

\section*{CRediT authorship contribution statement}
{\bf J. H. Woo:} Methodology, Software, Investigation, Visualization, Writing – original draft, Writing - Review \& Editing.     
{\bf H. S. Lee:} Methodology, Formal analysis, Software, Investigation, Writing – original draft, Writing - Review \& Editing.
{\bf J.-Y. Moon:} Conceptualization, Supervision, Methodology, Software, Investigation, Writing - Review \& Editing.
{\bf T.-W. Ko:} Conceptualization, Supervision, Methodology, Formal analysis, Investigation, Software, Visualization, Writing – original draft, Writing - Review \& Editing. 

\section*{Declaration of competing interest}
The authors declare that they have no known competing financial interests or personal relationships that could have appeared to influence the work reported in this paper.

\section*{Data availability}
No data was used for the research described in the article.

\section*{Declaration of generative AI and AI-assisted technologies in the writing process}
During the preparation of this work the authors used Claude.ai 3.5 Sonnet in order to improve language and readability. After using this tool/service, the authors reviewed and edited the content as needed and take full responsibility for the content of the publication.

\begin{acknowledgements}
	{This work was supported by the National Research Foundation of Korea (NRF) grant funded by the Korea government (MSIT) grant No.2019R1A2C2089463(HSL), the National Institute for Mathematical Sciences (NIMS) grant funded by the Korean government (No. NIMS-B24910000), and Institute for Basic Science (IBS) grant No. IBS-R015-Y3.
}
\end{acknowledgements}

\appendix
\section{\label{appendix:A} Stability Analysis of states using Ott-Antonsen (OA) approach}
Applying the Ott-Antonsen (OA) approach \cite{ott2008,ott2009} similarly to our model of Eq.~(\ref{eq:model}) as in Refs.~\cite{laing2009,wzou2020,cxu2021}, we derive an equation for the time evolution of the states of the system. We identify stationary states such as uniformly incoherent states and partially (fully) locked states, and conduct a linear stability analysis of these states.

\subsection{Finding stationary states using Ott-Antonsen (OA) ansatz}
In the continuum limit, the states of the system described by Eq.~(\ref{eq:model}) can be characterized by a probability density function (PDF), $f(\theta,K,t)$. For each value of $K$, the expression $f(\theta,K,t)d\theta$ represents the fraction of oscillators whose phases fall within the interval $[\theta, \theta+d\theta]$ at time $t$. The PDF $f(\theta,K,t)$ must satisfy the normalization condition:

\begin{equation}
\label{eq:normalization}
	\int_{0}^{2\pi} d\theta\, f(\theta,K,t) = 1 ~~{\rm for~all}~ K
\end{equation}
and the continuity equation:
\begin{eqnarray}
\label{eq:conti}
    \frac{\partial{f}}{\partial{t}} + \frac{\partial}{\partial{\theta}}(fv) = 0
\end{eqnarray}
with $v$ given by
\begin{eqnarray}
\label{eq:v}
v &=& \omega + Kc_0 + K \int_0^\infty dK' \int_0^{2\pi} d\theta' \sin(\theta'-\theta-\beta) \nonumber \\
&& ~~~~~~~~~\times f(\theta',K',t) g(K'),
\end{eqnarray}
where $g(K)$ is the PDF for coupling strength $K$.

Using a complex order parameter $z$ defined by
\begin{eqnarray}
z(t) &=& r(t){\rm e}^{i\Theta(t)} \nonumber \\
     &=&  \int_0^\infty dK \int_0^{2\pi} d\theta' {\rm e}^{i\theta'} f(\theta',K,t) g(K),
\label{eq:z}
\end{eqnarray}
we can rewrite Eq.~(\ref{eq:v}) as
\begin{eqnarray}
v &=& \omega + Kc_0 + \frac{K}{2i}\big [ z {\rm e}^{-i(\theta+\beta)} - z^* {\rm e}^{i(\theta+\beta)} \big].
\end{eqnarray}
The PDF $f(\theta, K, t)$ can be written as 
\begin{eqnarray}
f(\theta, K, t) = \frac{1}{2\pi}\left \{ 1+ \left [ \sum_{n=1}^\infty a_n(K,t){\rm e}^{in\theta}+ c.c. \right ] \right \},~~~~
\end{eqnarray}
where $c.c.$ stands for the complex conjugate of the preceding term.

With OA ansatz \cite{ott2008} $a_n(K,t) = a(K,t)^n$, the PDF becomes
\begin{eqnarray}
f(\theta, K, t) = \frac{1}{2\pi}\left \{ 1+ \left [ \sum_{n=1}^\infty a(K,t)^n{\rm e}^{i n\theta} + c.c. \right ] \right \}.~~~~
\label{eq:f_with_oa}
\end{eqnarray}

To guarantee that the amplitude of each mode does not diverge, we assume that $|a(K,t)| \le 1$.

Applying Eq. (\ref{eq:f_with_oa}) to the continuity equation in Eq.~\eqref{eq:conti}, we get an equation describing the time evolution of $a(K,t)$ as follows. 
\begin{eqnarray}
    \frac{\partial a}{\partial t} &=&  - i( \omega + Kc_0)a \nonumber \\ 
	&&+ \frac{K}{2}\biggl [ r{\rm e}^{i(\beta-\Theta)} - r{\rm e}^{-i(\beta-\Theta)} a^2 \biggr ],
\label{eq:a}
\end{eqnarray}
where Eq. (\ref{eq:z}) for the definition of the complex order parameter is used.

Inserting $f(\theta,K,t)$ of Eq.~(\ref{eq:f_with_oa}) into Eq.~(\ref{eq:z}), we obtain an expression for the complex order parameter with $a(K,t)$ : 
\begin{eqnarray}
z(t) &=& r(t){\rm e}^{i\Theta(t)} \nonumber \\
     &=& \int_0^\infty dK ~a^*(K,t) g(K).
\label{eq:z_with_a}
\end{eqnarray}

We aim to identify stationary states where $r(t)$ remains constant over time, while $\Theta(t)$ exhibits a linear time dependence: $\Theta(t) = \Phi_0 + \Omega t$. Here, $\Phi_0$ represents an arbitrary constant offset, which we can conveniently set to zero without affecting the generality of our analysis.

To achieve this, we introduce substitutions. We define  $a(K,t)=\tilde{a}(K,t) {\rm e}^{-i\Omega t}$ and $\Theta(t) = \Phi(t) + \Omega t$ with $\Phi(0) = 0$. These modified expressions will then be substituted into Eqs. (\ref{eq:a}) and (\ref{eq:z_with_a}), allowing us to solve for the desired stationary states.
This gives
\begin{eqnarray}
    \frac{\partial \tilde{a}}{\partial t} &=&  - i( \Delta + Kc_0)\tilde{a} \nonumber \\
    &&+ \frac{K}{2}\biggl [ r{\rm e}^{i(\beta-\Phi)}
    - r{\rm e}^{-i(\beta-\Phi)} {\tilde{a}}^2 \biggr ],
\label{eq:tilde_a}
\end{eqnarray}
where $\Delta = \omega - \Omega$, and 
\begin{eqnarray}
r(t){\rm e}^{i\Phi(t)} = \int_0^\infty dK ~\tilde{a}^*(K,t) g(K). 
\label{eq:R_with_tilde_a}
\end{eqnarray}

For the stationary states, $\tilde{a}$, $r$, and $\Phi$ must remain constants over time. To satisfy the condition of time-invariance for $\tilde{a}$, the right-hand side of Eq.~(\ref{eq:tilde_a}) must equate to zero. We obtain the following:
\begin{eqnarray}
    Kr{\rm e}^{-i\beta} \tilde{a}^2 
	+ 2iZ(K)\tilde{a} 
    - Kr{\rm e}^{i\beta} = 0,
\label{eq:steady_condition}
\end{eqnarray}
where $Z(K) = \Delta + Kc_0$ and $\Phi$ is set to zero. This equation has two solutions for $\tilde{a}$ with a given value of $K$:
\begin{eqnarray}
\label{eq:tilde_a_pm}
\tilde{a} &=& \tilde{a}_{\pm} \\
&=&\left\{
    \begin {aligned}
	 & \frac{-iZ(K)\pm \sqrt{K^2r^2-Z(K)^2}}{Kr{\rm e}^{-i\beta}} &{\rm for~} K \in \mathcal{D}_l,\nonumber \\
	 &  \frac{-iZ(K)\pm i\sqrt{Z(K)^2-K^2r^2}}{Kr{\rm e}^{-i\beta}}&{\rm for~} K \in \mathcal{D}_d,                  
    \end{aligned}
\right.
\end{eqnarray}
where $\mathcal{D}_{l} = \{K : K r > |Z(K)|\,\}$ and $\mathcal{D}_{d} = \{K : K r < |Z(K)|\,\}$.
The condition $|a(K,t)|\le 1$ requires the sign in front of the square root should be ${\rm sign}(Z(K))$ for $K\in \mathcal{D}_d$. In contrast, both of the sign in front of the square root can be possible for  $K\in \mathcal{D}_l$, since $|\tilde{a}_\pm| = 1$.

We get the following from Eq. (\ref{eq:R_with_tilde_a}):
\begin{eqnarray}
\label{R2_with_pm}
    r^2 {\rm e}^{i\beta}
    &=& i\int_{\mathcal{D}_{tot}} dK\: \frac{g(K)Z(K)}{K} \nonumber \\
    &&\pm \int_{\mathcal{D}_{l}} dK\: \frac{g(K)\sqrt{K^2r^2 - Z(K)^2}}{K} \\
    &&-i\int_{\mathcal{D}_{d}} dK\: \frac{g(K){\rm sign}(Z(K))\sqrt{Z(K)^2-K^2r^2}}{K},\nonumber
\end{eqnarray}
where $\mathcal{D}_{tot}$ is the total range of $K$. With $\beta \in [0,\pi/2)$, the real part of the left-hand side of Eq.~(\ref{R2_with_pm}) is nonnegative and thus only the positive branch of $\tilde{a}$ for $K \in \mathcal{D}_l$ is possible. Thus,
\begin{eqnarray}
\label{eq:tilde_a_final}
\tilde{a} &=& 
\left\{
    \begin {aligned}
	 & \frac{-iZ(K)+\sqrt{K^2r^2-Z(K)^2}}{Kr{\rm e}^{-i\beta}}~~{\rm for~} K \in \mathcal{D}_l,~~~~ \\
	 & \frac{-iZ(K)+i {\rm sign}(Z(K))\sqrt{Z(K)^2-K^2r^2}}{Kr{\rm e}^{-i\beta}}~~~~\\
	 &~~~~~~~~~~~~~~~~~~~~~~~~~~~~~~~~~~~~~~~~{\rm for~} K \in \mathcal{D}_d,               
    \end{aligned}
\right.
\end{eqnarray}

Using this, we derive the same self-consistency equation from Eq.~(\ref{eq:R_with_tilde_a}) as Eq.~(\ref{eq:self_consistency}).
\begin{eqnarray}
	r^2 {\rm e}^{i\beta} &=& i\int_{\mathcal{D}_{tot}} dK \: \frac{g(K)Z(K)}{K}  
\nonumber \\
	&&+\int_{\mathcal{D}_l} dK\: \frac{g(K)\sqrt{K^2 r^2 -Z(K)^2}}{K} \\
	&&-i\int_{\mathcal{D}_d} dK\: \frac{g(K){\rm sign}\big(Z(K)\big)\sqrt{Z(K)^2-K^2 r^2}}{K}. \nonumber
\label{eq:self_consistency_oa}
\end{eqnarray}
We can numerically determine $r$ and $\Delta$ from the above equation, and thus $\tilde{a}$ from Eq.~(\ref{eq:tilde_a_final}), which specifies the stationary states.

\subsection{Stability analysis of uniformly incoherent state}
First, we perform the linear stability analysis of the uniformly incoherent state, which is described by $f(\theta,K,t) = \frac{1}{2\pi}$, corresponding to $\tilde{a}(K,t) = 0$. This state is characterized by $r=0$ and $\Delta=0$.

Using a perturbation to $\tilde{a}(K,t) = 0$,
\begin{equation}
    \tilde{a}(K,t) = 0 + \varepsilon \eta(K,t),
\end{equation}
where $0<\varepsilon \ll 1$, we obtain
\begin{equation}
r(t) {\rm e}^{i\Phi(t)} = \varepsilon \int_{0}^{\infty} dK \, \eta^*(K,t) g(K)
\end{equation}
from Eq.~(\ref{eq:R_with_tilde_a}).

With the introduction of an operator $\hat{g}x \equiv \int_{0}^{\infty} dK g(K)x(K,t)$ and
the perturbation, we derive the following in the first order of $\varepsilon$ from Eq.~(\ref{eq:tilde_a}):
\begin{align}
    \frac{\partial \eta}{\partial t} &= -iKc_0 \eta + \frac{K{\rm e}^{i\beta}\hat{g}\eta}{2}.
\end{align}
Let $\hat{A}$ denote the operator defined by 
\begin{align}
\hat{A}x \equiv -iKc_0 x + \frac{K{\rm e}^{i\beta}\hat{g}x}{2}.
\label{eq:hat_A}
\end{align}

An eigenvalue $\lambda$ of the operator $\hat{A}$ and the corresponding eigenvector $v$ satisfy the following eigenvalue equation 
\begin{eqnarray}
	\hat{A} v = \lambda v
\label{eq:A_eigenvalue_eq}
\end{eqnarray}
and $\lambda$ determines the stability of the uniformly incoherent state.

From Eqs.~(\ref{eq:hat_A}) and (\ref{eq:A_eigenvalue_eq}), 
\begin{align}
-iKc_0 v + \frac{K{\rm e}^{i\beta}\hat{g}v}{2} = \lambda v
\end{align}

The eigenvalue $\lambda = -iKc_0$ and the corresponding eigenvector $v$ satisfying $\hat{g}v=0$ do not contribute to the instability of the uniformly incoherent state, since ${\rm Re}(\lambda) = 0$.

For $\lambda \neq -iKc_0$, which implies $\hat{g}v \neq 0$, 
\begin{eqnarray}
	v = \frac{K{\rm e}^{i\beta}\hat{g}v}{2(\lambda+iKc_0)}.
\label{eq:v_gv_A}
\end{eqnarray}

Applying the operator $\hat{g}$ to both sides of Eq.~(\ref{eq:v_gv_A}) yields
\begin{eqnarray}
    \hat{g}v = \frac{{\rm e}^{i\beta}\hat{g}v}{2}\int_0^\infty dK \,\frac{Kg(K)}{\lambda+iKc_0}
\end{eqnarray}
Thus, we get
\begin{eqnarray}
	1 &=& \frac{{\rm e}^{i\beta}}{2} \int_0^\infty dK\,\frac{Kg(K)}{\lambda+iKc_0}.
\end{eqnarray}

With $\lambda = \mu - i\nu$, we obtain the same equations for $\mu$ and $\nu$ as in Eqs.~(\ref{eq:stability_a}) and (\ref{eq:stability_b}).   
\begin{subequations}
\begin{eqnarray}
\cos \beta &=& \frac{1}{2} \int_0^\infty dK \, \frac{\mu K g(K)}{\mu^2+(Kc_0-\nu)^2}, \\
\sin \beta &=& \frac{1}{2} \int_0^\infty dK \, \frac{K g(K)(Kc_0-\nu)}{\mu^2+(Kc_0-\nu)^2}. 
\end{eqnarray}
\end{subequations}

Following the same approach using Eqs.~(\ref{eq:stability_a}) and (\ref{eq:stability_b}) as in subsection, we can identify the critical $c_0$ above which uniformly incoherent states is obtained.

\subsection{Stability analysis of partially (fully) locked states}
Next, we perform the linear stability analysis of a partially (fully) locked state, which is described by $\tilde{a}(K,t) = \tilde{a}_0(K)$. This state is characterized by $r=r_0(>0)$, $\Theta(t) = \Omega t$ and $\Delta=\Delta_0$. $\tilde{a}_0(K)$, $r_0$, and $\Delta_0$ satisfy Eq.~(\ref{eq:steady_condition}).

Using a perturbation to $\tilde{a}(K,t) = \tilde{a}_0(K)$,
\begin{equation}
    \tilde{a}(K,t) = \tilde{a}_0(K) + \varepsilon \eta(K,t),
\label{eq:tilde_a0_plus_epsilon}
\end{equation}
where $0<\varepsilon \ll 1$, we obtain

\begin{equation}
r(t){\rm e}^{i\Phi(t)} = r_0(t) + \varepsilon \hat{g}\eta^* 
\label{eq:R0_plus_epsilon}
\end{equation}
from Eq.~(\ref{eq:R_with_tilde_a}).

Using Eqs.~(\ref{eq:tilde_a0_plus_epsilon}) and (\ref{eq:R0_plus_epsilon}), 
we derive the following evolution equation for $\eta$ in the first order of $\varepsilon$ from Eq.~(\ref{eq:tilde_a}):
\begin{eqnarray}
    \frac{\partial \eta}{\partial t} &=&  
    -iZ_0(K)\eta + \frac{K}{2} \biggl [{\rm e}^{i\beta}\hat{g}\eta \nonumber \\
    &&~~~~~~~~~~~~ - {\rm e}^{-i\beta}(\tilde{a}_0^2 \hat{g}\eta^* + 2 \tilde{a}_0 r_0 \eta )\biggr ], 
	\label{eq:eta_eq_partial_locked}
\end{eqnarray}
where $Z_0(K)=\Delta_0+Kc_0$.

From Eq.~(\ref{eq:eta_eq_partial_locked}), we get  
\begin{eqnarray}
    \frac{\partial V}{\partial t} &=& \mathbf{M}V +\mathbf{Q}\hat{g}V, %\nonumber \\
    %&\equiv & \mathcal{L} V
\end{eqnarray}
where $V=(\eta, \eta^*)^T$ and
\begin{eqnarray}
	\mathbf{M} &=& 
\begin{pmatrix}
m_{11}(K) & 0 \\
0 & m_{11}^*(K)
\end{pmatrix}, \\
	\mathbf{Q} &=& 
\begin{pmatrix}
q_{11}(K) & q_{12}(K) \\
	q_{12}^*(K) & q_{11}^*(K)
\end{pmatrix},
\end{eqnarray}
with 
\begin{eqnarray}
\label{eq:m11}
m_{11}(K) &=& -iZ_0(K)-K{\rm e}^{-i\beta}\tilde{a}_0 r_0 \\
&=&
\left\{
    \begin {aligned}
         & -\sqrt{K^2r_0^2-Z_0(K)^2}~~{\rm for~}K \in \mathcal{D}_l, \nonumber \\
         &  {-i {\rm sign}(Z_0(K))\sqrt{Z_0(K)^2-K^2r_0^2}} \nonumber \\
	 &~~~~~~~~~~~~~~~~~~~~~~~~~~~~~~{\rm for~}K \in \mathcal{D}_d, \nonumber \\                
    \end{aligned}
\right.\\
q_{11}(K) &=& K{\rm e}^{i\beta}/2,\\
q_{12}(K) &=& -K{\rm e}^{-i\beta}\tilde{a}_0^2/2.
\end{eqnarray}

Let $\hat{B}$ denote the operator defined by 
\begin{align}
	\hat{B}X \equiv \mathbf{M} X + \mathbf{Q}\hat{g}X.
\label{eq:hat_B}
\end{align}
with $X=(x,x^*)^T$.

An eigenvalue $\lambda$ of the operator $\hat{B}$ and the corresponding eigenvector $V=(v,v^*)^T$ satisfy the following eigenvalue equation
\begin{eqnarray}
    \hat{B} V = \lambda V
\label{eq:B_eigenvalue_eq}
\end{eqnarray}
and $\lambda$ determines the stability of the partially(fully) locked state.

The eigenvalues $\lambda = m_{11}(K)$ and $\lambda = m_{11}^*(K)$, and the corresponding eigenvectors $V$ satisfying $\hat{g}V=0$ do not contribute to the instability of the partially(fully) locked state, because ${\rm Re}(\lambda) = {\rm Re}(m_{11}(K))$ is either negative or zero according to Eq.~(\ref{eq:m11}).

Thus, let us focus on the eigenvalues with the conditions $\lambda \neq m_{11}(K)$ and $\lambda \neq m_{11}^*(K)$. Since $\lambda \mathbf{I} - \mathbf{M}$ is invertible with such eigenvalues, we can obtain 
\begin{eqnarray}
    V = (\lambda \mathbf{I} - \mathbf{M})^{-1}\mathbf{Q}\hat{g}V
\label{eq:V_gV_B}
\end{eqnarray}
from Eqs.~(\ref{eq:hat_B}) and (\ref{eq:B_eigenvalue_eq}).

Applying the operator $\hat{g}$ to both sides of Eq.~(\ref{eq:V_gV_B}), we obtain
\begin{eqnarray}
    (\mathbf{I} - \mathbf{J})\hat{g}V = 0,
\end{eqnarray}
where
\begin{eqnarray}
	\mathbf{J} &=& \hat{g}(\lambda \mathbf{I} - \mathbf{M})^{-1}\mathbf{Q} \nonumber \\
&=& \hat{g}
\begin{pmatrix}
[\lambda-m_{11}]^{-1} & 0 \\
0 & [\lambda-m_{11}^*]^{-1}
\end{pmatrix}
\begin{pmatrix}
q_{11} & q_{12} \\
q_{12}^* & q_{11}^*
\end{pmatrix} \nonumber \\
&=& 
\begin{pmatrix}
J_{11}(\lambda) & J_{12}(\lambda) \\
J_{21}(\lambda) & J_{22}(\lambda)
\end{pmatrix}
\end{eqnarray}
with
\begin{eqnarray}
	\label{eq:J11}
	J_{11}(\lambda) &=& \int_0^\infty dK \: \frac{q_{11}(K)g(K)}{\lambda-m_{11}(K)}, \\
	\label{eq:J12}
	J_{12}(\lambda) &=& \int_0^\infty dK \: \frac{q_{12}(K)g(K)}{\lambda-m_{11}(K)}, \\
	\label{eq:J21}
	J_{21}(\lambda) &=& \int_0^\infty dK \: \frac{q_{12}^*(K)g(K)}{\lambda-m_{11}^*(K)}, \\
	\label{eq:J22}
	J_{22}(\lambda) &=& \int_0^\infty dK \: \frac{q_{11}^*(K)g(K)}{\lambda-m_{11}^*(K)}. 
\end{eqnarray}

The characteristic equation of the operator $\hat{B}$ becomes  
\begin{eqnarray}
    {\rm det} \biggl [\mathbf{I}-\mathbf{J}(\lambda) \biggr] = 0,
\label{eq:char_eq_B}
\end{eqnarray}
which is equivalent to 
\begin{eqnarray}
	(1-J_{11}(\lambda))(1-J_{22}(\lambda)) - J_{12}(\lambda)J_{21}(\lambda) = 0.~~
\label{eq:char_eq_B_equiv}
\end{eqnarray}
The partially (fully) locked state is unstable if there is an eigenvalue $\lambda$ satisfying Eq.~(\ref{eq:char_eq_B_equiv}) and the condition $\mathrm{Re}(\lambda) > 0$.  We numerically check the existence of such an eigenvalue and determine the stability by searching for it in the right half-plane of the complex plane where $\mathrm{Re}(\lambda)>0$. During this search, there are no singularity-related problems since the integrands of Eqs.~(\ref{eq:J11}) - (\ref{eq:J22}) do not have singularities and are continuous for $\mathrm{Re}(\lambda)>0$, ensuring the existence of the integrations. 
If such an eigenvalue is not found, the partially locked state and the fully locked state are linearly neutrally stable and linearly stable, respectively.

\end{document}